\documentclass[12pt]{amsart}
\usepackage{amsfonts,amsmath,amssymb,amsthm}
\usepackage[all]{xy}
\usepackage{anysize}
\usepackage{graphicx}

\marginsize{1in}{1in}{1in}{1in}

\def\Q{\mathbf{Q}}
\def\Z{\mathbf{Z}}
\def\F{\mathbf{F}}
\def\H{\mathbf{H}}
\def\R{\mathbf{R}}
\def\C{\mathbf{C}}
\def\P{\mathbf{P}}
\def\Fp{{\F}_p}
\def\Fpbar{{\overline\F_p}}
\def\Fl{{\F}_l}
\def\A{\mathbf{A}}
\def\gp{{\mathfrak p}}
\def\gP{{\mathfrak P}}
\def\Otilde{\widetilde O}

\def\Gal{\text{\rm Gal}}
\def\Jac{\text{\rm Jac}}
\def\disc{\text{\rm disc}}
\def\Stab{\text{\rm Stab}}
\def\Sp{\text{\rm Sp}}
\def\Sym{\text{\rm Sym}}
\def\Cl{\text{\rm Cl}}

\def\Im{\text{\rm Im}}

\def\Tr{\text{\rm Tr}}
\def\Spec{\text{\rm Spec}}

\def\gC{\mathfrak{C}}
\def\sA{\mathfrak{A}}
\def\endproof{\hfill$\square$\par}
\def\isar{\ \smash{\mathop{\longrightarrow}\limits^{\thicksim}}\ }
\gdef\mapright#1{\ \smash{\mathop{\longrightarrow}\limits^{#1}}\ }

\newcommand{\ZZ}{\mathbf{Z}}
\newcommand{\QQ}{\mathbf{Q}}

\newcommand{\M}{\mathcal{M}}
\newcommand{\Am}{\mathcal{A}} 
\newcommand{\Order}{\mathcal{O}}

\newcommand{\smat}[4]{\left( \smallmatrix #1 & #2 \\ #3 & #4
    \endsmallmatrix \right)}

\newcommand{\Hom}{{\rm Hom}}
\newcommand{\End}{{\rm End}}

\def\ga{\mathfrak{a}}
\def\gb{\mathfrak{b}}

\def\D{{\mathfrak{D}}}
\def\Ibar{{\overline{I}}}

\newtheorem{theorem}{Theorem}[section]
\newtheorem{lemma}[theorem]{Lemma}

\newcommand{\Remark}{\noindent{\bf Remark.}\ }
\newcommand{\Proof}{\noindent{\bf Proof.}\ }
\newtheorem{algorithm}[theorem]{Algorithm}
\newtheorem{corollary}[theorem]{Corollary}
\newtheorem{heuristic}[theorem]{Heuristic}

\theoremstyle{definition}
\newtheorem{remark}[theorem]{Remark}
\newtheorem{example}[theorem]{Example}

\numberwithin{equation}{section}

\begin{document}

\title{Explicit CM-theory for level 2-structures on abelian surfaces}
\author{Reinier Br\"oker, David Gruenewald, and Kristin Lauter}
\address{Brown University, Department of Mathematics, Providence, RI 02912, USA}
\email{reinier@math.brown.edu}
\address{Laboratoire de Math\'ematiques Nicolas Oresme, Universit\'e �de Caen Basse-Normandie, 14032 Caen cedex, France}
\email{davidg@maths.usyd.edu.au}
\thanks{The second author thanks Microsoft Research, where this research was
undertaken, for its hospitality.}
\address{Microsoft Research, One Microsoft Way, Redmond, WA 98052, USA}
\email{klauter@microsoft.com}


\begin{abstract}
For a complex abelian surface $A$ with endomorphism ring isomorphic to the
maximal order in a quartic CM-field $K$, the Igusa invariants $j_1(A),
j_2(A),j_3(A)$ generate an unramified abelian extension of the reflex field of~$K$.  In
this paper we give an explicit geometric description of the Galois action of 
the class group of this reflex field on $j_1(A),j_2(A),j_3(A)$. Our
description can be expressed by maps between various Siegel modular
varieties, and we can explicitly compute the action for ideals of small norm. 
We use the Galois action to modify the CRT method for computing Igusa class 
polynomials, and our run time analysis shows that this yields a significant
improvement.
Furthermore, we find cycles in isogeny graphs for abelian surfaces, thereby
implying that the `isogeny volcano' algorithm to compute endomorphism rings of
ordinary elliptic curves over finite fields does not have a straightforward
generalization to computing endomorphism rings of abelian surfaces over finite
fields.
\end{abstract}

\maketitle

\section{Introduction} \label{sec1}
Class field theory describes the abelian extensions of a given number 
field~$K$. For $K=\Q$, the Kronecker-Weber theorem tells us that every
abelian extension of $K$ is contained in a {\it cyclotomic extension\/}. In
1900, Hilbert asked for a similar `explicit description' for higher degree
number fields. This is known as Hilbert's 12th problem, and it is still largely unsolved. 

Besides $K=\Q$, the answer is only completely known for imaginary quadratic
fields. In this case, the solution is provided by {\it complex
multiplication\/} theory, see e.g.~\cite[Ch.~2]{Sil}. The techniques used can be
generalized to {\it CM-fields\/}, i.e., imaginary quadratic extensions of
totally real fields. However, for general CM-fields we do not always get an
explicit description of the {\it maximal\/} abelian extension. From a
computational perspective, the case of general CM-fields is far less developed
than the imaginary quadratic case.

In this article, we solely focus on degree $4$ primitive CM-fields~$K$. For such
fields, invariants of principally polarized abelian surfaces (p.p.a.s.) with
endomorphism ring isomorphic to the maximal order $\Order_K$ of $K$ generate a
subfield of the Hilbert class field of the {\it reflex field\/} of~$K$ (a
degree 4 subfield of the normal closure of~$K$). 
To explicitly compute the
resulting extension, we compute an {\it Igusa class polynomial\/} 
$$ 
P_K = \prod_{\{A\ p.p.a.s. \mid \End(A) = \Order_K\}/\cong} \left(X - j_1(A)\right)
\in\Q[X].  $$ 
Here, $j_1$ is one of the {\it three\/} Igusa invariants of~$A$.
A contrast with the case of imaginary quadratic fields -- where we compute the
{\it Hilbert class polynomial\/} -- is that the polynomial $P_K$ has rational
coefficients which are not integers in general, and it need not be irreducible
over~$\Q$.

There are three methods to explicitly compute the polynomial $P_K$: complex 
analytic evaluation of the invariants~\cite{Spallek, van Wamelen, Weng}, the CRT method 
using finite field arithmetic~\cite{EisentragerLauter} and the 
computation of a canonical lift~\cite{2adicFiveAuthor, 3adic} using $p$-adic 
arithmetic for $p=2,3$. However, none of these three approaches 
exploit the Galois action of the maximal abelian extension of the reflex field 
on the set of 
principally polarized abelian surfaces with 
endomorphism ring~$\Order_K$. The goal of this article is to make this Galois
action explicit and give a method to compute it. 

Our algorithm to compute
the Galois action significantly speeds up the CRT-approach described in~\cite{EisentragerLauter} to compute Igusa class polynomials and it can be
used to improve the $3$-adic approach~\cite{3adic} as well. The improvement
in computing Igusa class polynomials parallels the improvements given 
in~\cite{BBEL} for computing Hilbert class polynomials. Our run time analysis
is similar to the analysis in~\cite{BBEL}. Contrary to the genus $1$ algorithm
however, the genus $2$ algorithm is not quasi-linear in the size of the
output. We suggest further refinements that might yield a quasi-linear
algorithm as area of further study in Section~6.

Besides speeding up the computation of Igusa class polynomials, our algorithm
gives a method of computing {\it isogenous\/} abelian surfaces
over finite fields. Computing an isogeny is a basic computational problem
in arithmetic geometry, and we expect that our algorithm can be used in a 
variety of contexts, ranging from point counting on Jacobians of curves to
cryptographic protocols. 

Our computations naturally lead to the study of the $(l,l)$-isogeny graph
of abelian surfaces over finite fields. For ordinary elliptic curves, the $l$-isogeny
graph looks like a `volcano' and this observation forms the heart of the
algorithm~\cite{KohelThesis} to compute the endomorphism ring of an ordinary elliptic curve over a finite field. We show that for abelian surfaces, the $(l,l)$-isogeny graph 
does {\it not\/} have a volcano shape. This shows that a straightforward
generalization of the elliptic curve algorithm to abelian surfaces does not work. 

The structure of this paper is as follows. In Section \ref{sec2} we recall the
basic facts of complex multiplication theory and background on CM abelian surfaces and their invariants. In Section \ref{sec3} we describe the 
Galois action on the set of isomorphism classes of abelian surfaces with CM 
by $\Order_K$ in a `geometric way'. 
Our algorithm to compute this action is intrinsically linked to Siegel 
modular functions of higher level.
Section \ref{sec4} gives the definitions and properties of the four Siegel 
modular functions that we use.
The algorithm to compute the Galois action is detailed 
in Section~\ref{sec5} and we apply it in Section~\ref{sec6} where we improve 
the method to compute an Igusa class polynomial modulo a prime~$p$. We give
a detailed run time analysis of our algorithm in Section~6 as well.
We illustrate our approach with various detailed examples in 
Section~\ref{sec7}. A final Section~\ref{sec8} contains the obstruction to the 
volcano picture for abelian surfaces.

\section{CM abelian surfaces} \label{sec2}

\subsection{CM Theory}\par
\ \par
In this section we recall the basic facts of CM-theory for higher dimensional
abelian varieties. Most of the material presented in this section is an
adaptation to our needs of the definitions and proofs of Shimura's~\cite{Shim}
and Lang's~\cite{LaCM} textbooks.

We fix an embedding of $\overline\Q \hookrightarrow \C$. By a real number
field, we mean a field that is fixed by complex conjugation. With this
convention, a {\it CM field\/} $K$ is a totally imaginary quadratic extension of a totally 
real number field. Let $K^+$ denote the real quadratic subfield of $K$, and
let $n$ be the degree of $K^+$ over $\Q$. The $2n$ embeddings 
$K \hookrightarrow \overline\Q$ naturally come in pairs. Indeed, we can choose $n$
embeddings $\Phi = \{\varphi_1,\ldots,\varphi_{n}\}$ such that we have
$\Hom(K,\overline\Q) = \Phi \cup \overline\Phi$. We call such a set $\Phi$ a 
{\it CM-type\/} for~$K$, and we interpret a CM-type in the natural way as
a map $K \hookrightarrow \C^n$. 

If $\Phi$ {\it cannot\/} be obtained as a lift of a CM-type of a CM-subfield
of $K$, then we call $\Phi$ {\it primitive\/}. For instance, in the simplest
case $K^+=\Q$, CM-fields $K$ are imaginary quadratic and every choice
for $K \hookrightarrow \overline\Q$ determines a primitive CM-type. If $K$ has degree
four then every choice of a CM-type is primitive when $K$ does not contain
an imaginary quadratic field. It is not hard to show \cite[Sec.~8.4]{Shim} that this occurs exactly 
for $\Gal(L/\Q) = D_4, C_4$, where $L$ denotes the normal closure of~$K$. We
say that the field $K$ is primitive in this case. 

{\it In this article, we will only consider primitive quartic 
CM-fields~$K$.\/} For the remainder of this section, we fix such a field~$K$.
We say that a principally polarized abelian surface $A/\C$ has \emph{CM} by the
maximal order $\Order_K$ if there exists an isomorphism $\Order_K \isar \End(A)$.
The CM-type distinguishes these surfaces. More precisely, a surface $A$ that
has CM by $\Order_K$ has 
{\it type\/} $\Phi = \{\varphi_1,\varphi_2\}$ if the complex representation 
$R_\C$ of the endomorphism algebra $\End(A) \otimes_{\Z} \Q$ satisfies
$$
R_\C \cong \varphi_1 \oplus \varphi_2.
$$
One shows~\cite[Thm.\ 1.3.6]{LaCM} that a principally polarized abelian
surface that has CM by $\Order_K$ of type $\Phi$ is {\it simple\/}, i.e., is not
isogenous to the product of elliptic curves.

Let $\Phi$ be a CM-type for~$K$.  For an $\Order_K$-ideal $I$, the quotient
$A_I = \C^2/\Phi(I)$
 is an abelian surface of type $\Phi$ by~\cite[Thm.\ 4.1]{LaCM}. 
This surface
need not admit a principal polarization. The dual variety of $A_I$ is given
by $\widehat{A_I} = \C^2/\Phi(\Ibar^{-1} \D_K^{-1})$, where
$$
\D_K^{-1} = \{x \in K \mid \Tr_{K/\Q}(x\Order_K) \subseteq \Z \}
$$
is the inverse different and $\Ibar$ denotes the complex conjugate of~$I$.
If $\pi\in K$ satisfies $\Phi(\pi) \in (i\R_{>0})^2$ and
$\pi \D_K = (I \Ibar)^{-1} $,
then the map $A_I \rightarrow \widehat{A_I}$ given by $$
(z_1,z_2) \mapsto (\varphi_1(\pi) z_1, \varphi_2(\pi) z_2)
$$
is an isomorphism~(\cite[p.\ 102--104]{Shim}) and $A_I$ is principally
polarizable. All principally polarized abelian surfaces with CM by $\Order_K$ of
type $\Phi$ arise via this construction.

Let $L$ be the normal closure of $K$.
We extend $\Phi$ to a CM-type $\Phi'$ of $L$, and we define the {\it reflex
field\/}
$$
K_\Phi = \Q\Bigl(\left\{ \textstyle\sum_{\phi\in\Phi'} \phi(x) \mid x \in K\right\}\Bigr).
$$
The CM-type on $K$ induces a CM-type $f_\Phi = \{ \sigma^{-1}\!\!\!\mid_{
K_\Phi}\ :\, \sigma \in \Phi'\}$ of the reflex field~$K_\Phi$. The field
$K_\Phi$ is a subfield of $L$ of degree 4. In particular, it equals
$K$ in the
case $K$ is Galois. If $L/\Q$ is dihedral, then $K_\Phi$ and $K$ are not
isomorphic. However, the two different CM-types yield isomorphic reflex
fields in this case. Furthermore, we have
$$
(K_\Phi)_{f_\Phi} = K
$$
and the induced CM-type on $(K_\Phi)_{f_\Phi}$ equals $\Phi$.

An automorphism $\sigma$ of $K$ induces an isomorphism 
$(A,\Phi)\isar (A^\sigma,\Phi^\sigma)$ of CM abelian surfaces where 
$\Phi^\sigma=\{\varphi_1\sigma,\varphi_2\sigma\}$.  Thus two CM-types which 
are complex conjugates of each other produce the same sets of isomorphism 
classes of abelian surfaces. In the Galois case there is only one CM-type 
up to isomorphism and in the dihedral case there are two distinct CM-types.

\subsection{Igusa invariants}\par
\ \par
Any principally polarized abelian surface over $\C$ is of the form
$A_\tau = \C^2/(\Z^2+\Z^2\tau)$ where $\tau$ is an element of the {\it Siegel upper
half plane\/}
$$
\H_2 = \{ \tau \in \textrm{Mat}_{2}(\C) \mid  \tau \textrm{ symmetric, }  
\Im(\tau) \text{\rm\ positive definite} \}.
$$
The moduli space $\Am_2$ of principally polarized abelian surfaces 
is $3$-dimensional. We are mostly interested in the subspace $\M_2 \subset
\Am_2$ of Jacobians of curves. The structure of~$\M_2$ is well known, we
recall it for convenience. 
Let $Y^2 = a_6 X^6+\ldots+a_0 = f(X)$ be a genus $2$ curve and write $\alpha_1,\ldots, \alpha_6$ for the roots of~$f$.  For simplifying notation, let $(ij)$ denote the quantity $(\alpha_{k_i}-\alpha_{k_j})$ for a given ordering of the roots.
The {\it Igusa-Clebsch invariants\/} $I_2,I_4,I_6,I_{10}$ (denoted by $A,B,C,D$ in \cite[Sec.~3]{IgusaAVM}) are defined by 
\begin{align*}
I_2 &= a_6^{2} \sum_{15} (12)^2 (34)^2 (56)^2\\
I_4 &= a_6^4 \sum_{10} (12)^2 (23)^2 (31)^2 (45)^2 (56)^2 (64)^2\\
I_6 &= a_6^6 \sum_{60} (12)^2 (23)^2 (31)^2 (45)^2 (56)^2 (64)^2 (14)^2 (25)^2 (36)^2\\
I_{10} &= a_6^{10} \sum_{i<j} (ij) = a_6^{10} \disc(f),
\end{align*}
where we sum over all root orderings $\{\alpha_{k_i}\}$ which give distinct summands; the subscript
indicates the number of terms we sum over. 

\begin{theorem} The moduli space $\M_2$ is isomorphic to $$\{ [I_2 : I_4 : I_6 : I_{10}] 
\in \P_w^3(\C) \mid I_{10} \not = 0 \}$$ where $\P_w^3$ denotes weighted
projective space with weights $2,4,6$ and $10$. 
\end{theorem}
\Proof See \cite{IgusaAVM}. 
\endproof\par\ \par
\noindent
We note that the condition $I_{10} \not = 0$ ensures that the polynomial $f$
defining the genus 2 curve is separable. 

Instead of working with a subset of weighted projective space, many people
work with an affine subspace of $\M_2$. This non-weighted subspace is 
given by 
$$
(j_1,j_2,j_3) = \left(\frac{I_2^5}{I_{10}}, \frac{I_4 I_2^3}{I_{10}}, \frac{I_6I_2^2}{I_{10}}\right).
$$
The functions $j_i$ are commonly called the {\it Igusa functions\/}. We
remark that there are various definitions of these functions and there are
different conventions for which choice is the `right' one. Our functions
are the same as in e.g.~\cite{van Wamelen}. They
have the property that for $\tau,\tau'$ corresponding to Jacobians of curves,
the equality $j_i(\tau) = j_i(\tau') \not = 0$ for $i=1,2,3$ implies that $C$ and $C'$ are 
isomorphic. A detailed description on
computing $j_i(\tau)$ for a point $\tau\in\H_2$ can be found 
in~\cite{DupontThesis, Weng}.

A `weak version' of the main theorem of complex multiplication theory is
that, for a primitive quartic CM-field $K$, the Igusa invariants of an
abelian variety with CM by $\Order_K$ generate an unramified abelian extension
of a reflex field of~$K$. More precisely, we have the following result.

\begin{theorem}  \label{sec2:thm1}
Let $(K, \Phi)$ be a primitive quartic CM type. Let $I$ be an $\Order_K$-ideal with the property that there exists a principal polarization on $A_I = \C^2/\Phi(I)$. Then the 
field $K_\Phi(j_1(A_I),j_2(A_I), j_3(A_I))$ is a subfield of the Hilbert 
class field of $K_\Phi$. The polynomial
$$
P_K = \prod_{A}(X-j_1(A)),
$$
with $A$ ranging over the isomorphism classes of principally polarized abelian surfaces with 
endomorphism ring $\Order_K$, has rational coefficients. The same is true for the polynomials $Q_K,R_K$ 
for the $j_2$ and $j_3$-invariants.
\end{theorem}
\Proof
This is~\cite[Thm.~5.8]{Spallek}. \endproof\par
\ \par\noindent
We will see in Corollary~\ref{sec2:cardinality} that, for any primitive CM-type $\Phi$, there always exists an $\Order_K$-ideal $I$ such that $A_I$ is principally polarizable. 

\section{CM-action} \label{sec3}
Throughout this section, we let $K$ be a fixed primitive quartic CM-field. 
We also fix a CM-type $\Phi: K \hookrightarrow \C^2$ and let $A/\C$ be a 
principally polarized abelian surface that has complex multiplication by 
$\Order_K$ of CM-type $\Phi$. 
\subsection{Galois action of the class group} \label{sec3.1} \par
\ \par
We define a group $\gC(K)$ as
$$
\{ (\ga,\alpha) \mid \ga \hbox{\ a fractional\ } \Order_K\hbox{-ideal
with\ } \ga\overline\ga = (\alpha) \hbox{\ and\ } \alpha\in K^+ \hbox{\ 
totally positive} \}/\sim
$$
where two pairs $(\ga,\alpha)$ and $(\gb,\beta)$ are equivalent if and only if 
there exists a 
unit $u\in K^*$ with $\gb = u\ga$ and $\beta = u\overline u \alpha$. The
multiplication is defined componentwise, and $(\Order_K,1)$ is the neutral
element of $\gC(K)$.

The group $\gC(K)$ naturally acts on the finite set $S(K,\Phi)$ of isomorphism
classes of principally 
polarized abelian surfaces that have CM by $\Order_K$ of a given type $\Phi$. 
Indeed, any such surface is given by an ideal $I$ determining the variety 
and a `$\Phi$-positive' element $\pi\in K$ giving the principal polarization. 
We now put
$$
(\ga,\alpha) \cdot (I,\pi) = (\ga I, \alpha\pi)
$$
for $(\ga,\alpha) \in \gC(K)$. 
By~\cite[Sec.\ 14.6]{Shim}, the action of
$\gC(K)$ on $S(K,\Phi)$ is transitive and free. In particular, we have
$|\gC(K)| = |S(K,\Phi)|$.

The structure of the group $\gC(K)$ is best described by the following
theorem.
Denote by $\Cl^+(\Order_{K^+})$ the narrow class group of $\Order_{K^+}$ and write $(\Order_{K^+}^*)^+$ for the group of totally positive units of $\Order_{K^+}$.
\begin{theorem} \label{sec2:thm3}
Let $K$ be a primitive quartic CM-field. Then the sequence
$$
1 \mapright{} (\Order_{K^+}^*)^+/N_{K/K^+}(\Order_K^*) \mapright{u \mapsto (\Order_K,u)} 
\gC(K) \mapright{(\ga,\alpha) \mapsto \ga} \Cl(\Order_K) \mapright{N_{K/K^+}} 
\Cl^+(\Order_{K^+}) \mapright{} 1
$$
is exact.
\end{theorem}
\Proof The exactness at $(\Order_{K^+}^*)^+/N_{K/K^+}(\Order_K^*)$ is
\cite[Sec.~14.5]{Shim}. It
remains to show that the sequence is exact at $\Cl^+(\Order_{K^+})$. To prove this, we
first prove\footnote{We thank Everett Howe for suggesting this argument.} that
there is a {\it finite\/} prime that is ramified in $K/K^+$.

Suppose that $K/K^+$ is unramified at all finite primes. By genus theory,
we then have $K = K^+(\sqrt{n})$ with $n \in \Z$. However, $K$ then has
$\Q(\sqrt{n})$ as quadratic subfield and $K$ is a biquadratic field. This
contradicts our assumption that $K$ is primitive.

Because there is a finite prime of $K^+$ that ramifies in $K$, the extensions
$K/K^+$ and $H^+(K^+)/K^+$ are linearly disjoint. Here, $H^+$ denotes
the narrow Hilbert class field. By Galois theory, we then have
$$
\Gal(H(K)/K) \twoheadrightarrow \Gal(KH^+(K^+)/K) \isar \Gal(H^+(K^+)/
K^+)
$$
which proves the theorem. \endproof\par
\ \par
\Remark The surjectivity of the last arrow was also proved 
in~\cite[Lemma 2.1]{Kohel} under the assumption that there exists a finite
prime that ramifies in $K/K^+$. Our proof shows in fact that such a prime
always exists.

\noindent
\begin{corollary} \label{sec2:cardinality}
Let $K$ be a primitive quartic CM-field.
The set $S(K)$ of isomorphism classes of principally polarized abelian surfaces with CM by $\Order_K$ has cardinality
$$
|S(K)| = 
\begin{cases}
|\gC(K)| & \text{ if }\; \Gal(K/\Q)\cong C_4 \,, \\
2|\gC(K)| & \text{ if }\; \Gal(K/\Q)\cong D_4 \,. \\
\end{cases}
$$
\end{corollary}

\Proof
By Theorem~\ref{sec2:thm3}, the cardinality $|S(K,\Phi)| = |\gC(K)|$
is independent of the choice of a CM-type~$\Phi$.
If we let $n$ denote the number of CM-types up 
to conjugacy, then the theorem follows immediately from the 
equality $|S(K)|=n|S(K,\Phi)|$.
\endproof\par
\ \par
\noindent
The Galois group $\Gal(K_\Phi(j_1(A))/K_\Phi)$ acts
in the following way on the set $S(K,\Phi)$.  With $f_\Phi$ the
CM-type on $K_\Phi$ induced by $\Phi$, we define $N_\Phi: K_\Phi
\rightarrow K$ by 
$$
N_\Phi(x) = \prod_{\varphi \in f_\Phi} \varphi(x).
$$
For an $\Order_{K_\Phi}$-ideal $I$, the $\Order_K$-ideal $N_\Phi(I)$ is called the
{\it typenorm\/} of $I$. We get a natural map $m: \Cl(\Order_{K_\Phi}) \rightarrow
\gC(K)$ defined by
$$
m(\gp) = (N_\Phi(\gp), N_{K_\Phi/\Q}(\gp))
$$
for degree 1 prime representatives~$\gp$.
The Galois group of $K_\Phi(j_1(A))/K_\Phi$ is a quotient of 
$\Gal(H(K_\Phi)/K_\Phi) \cong \Cl(\Order_{K_\Phi})$, and by~\cite[Sec.\ 15.2]{Shim}
the induced map 
$$
m : \Gal(K_\Phi(j_1(A))/K_\Phi) \rightarrow \gC(K)
$$
is {\it injective\/}. This describes the Galois action. Indeed, the group $\gC(K)$ acts on the set of all principally polarized abelian surfaces that have CM by $\Order_K$, and $m$ maps the Galois group injectively into $\gC(K)$. 
In Example~\ref{sec7:ex2} we will
see that the natural map $\Cl(\Order_{K_\Phi}) \rightarrow \gC(K)$ need not be 
injective. 

The typenorm
 can be defined in a slightly different way as well. If $K/\Q$ is Galois
with $\Gal(K/\Q) = \langle \sigma \rangle$, then for $\Phi = \{1,\sigma\}$ 
we have $N_\Phi(\gp) = 
\gp^{1+\sigma^3}$. If $K$ is not Galois, then we have $N_\Phi(\gp) = 
N_{L/K}(\gp\Order_L)$. We will use this description both for actual computations
and in the proof of Lemma~\ref{sec6:lemgenus}.

In the remainder of Section~\ref{sec3.1}, we provide the theoretical framework that
will allow us to explicitly compute the CM-action.
Let $I$ be a $\Order_{K_\Phi}$-ideal of norm $l$. We assume for simplicity that
$l$ is prime. We have $m(I) = (N_\Phi(I),l) = (J, l) \in \gC(K)$, 
where $J$ is an $\Order_K$-ideal of norm~$l^2$.
\begin{lemma} \label{sec3:lem1}
Let $I$ be an $\Order_{K_\Phi}$-ideal of prime norm $l$ with typenorm 
$N_\Phi(I) = J \subset \Order_K$. Then $J$ divides $(l) \subset
\Order_K$.
\end{lemma}
\Proof This follows from the relation $N_\Phi(I) \overline{N_\Phi(I)} = 
N_{K_\Phi/\Q}(I) = l$.
\endproof\par
\ \par 
\noindent
For an $\Order_K$-ideal $M$, we define the `$M$-torsion' of the abelian surface
$A$ by
\begin{equation*}
A[M] = \{ P \in A(\C) \mid \forall \alpha\in M: \alpha(P) = 0 \}.
\end{equation*}
We assume here that we have {\it fixed\/} an isomorphism $\End(A) \isar \Order_K$,
meaning that $M$ is an $\End(A)$-ideal as well. If $M$ is 
generated by an integer $n$, then $A[M]$ equals the $n$-torsion $A[n]$. 

Lemma~\ref{sec3:lem1} implies that $A[J]$ is a 2-dimensional subspace of the
$l$-torsion $A[l]$ of~$A$. The polarization of $A$ induces 
a symplectic form on $A[l]$, and $A[l]$ is a {\it symplectic\/} vector space 
of dimension 4 over the finite field~$\Fl$. By CM-theory we know that the quotient 
$$
A / A[J]
$$
is again a {\it principally polarized\/} abelian surface. By~\cite[Sec.\ 23]{Mum}, this implies that $A[J]$ is an {\it isotropic\/} 
$2$-dimensional subspace of~$A[l]$, i.e., the symplectic form vanishes 
on $A[J]$. We recall that an isogeny $A \rightarrow B$ between principally 
polarized abelian surfaces whose kernel is a 
$2$-dimensional isotropic subspace of $A[l]$ is called an {\it $(l,l)$-isogeny\/},
and 
$$
A \rightarrow A/A[J]
$$
is an example of an $(l,l)$-isogeny. 

The moduli space of all pairs $(A,G)$, with $A$ a principally polarized
abelian surface over $\C$ and $G$ a $2$-dimensional isotropic subspace of $A[l]$
can be described by an ideal $V(l) \subset \Q[X_1,Y_1,Z_1,X_2,Y_2,Z_2]$. More
precisely, the variety corresponding to $V(l)$ equals the Siegel modular
variety $Y_0^{(2)}(l)$ studied e.g.\ in~\cite{BrokerLauterPrePub}. As a complex Riemann 
surface, we have 
$$
Y_0^{(2)}(l) = \Gamma_0^{(2)}(l) \backslash \H_2,
$$
with
$$ 
\Gamma_0^{(2)}(l) = \left\{ \smat{a}{b}{c}{d} \in
\Sp_4(\Z) \mid c \equiv 0_2 \bmod l \right\}.
$$

If we specialize $V(l)$ at a point $(X_1,Y_1,Z_1) = (j_1(A),j_2(A),j_3(A))$ then
the resulting ideal $V'(l)$ is $0$-dimensional. The corresponding variety is a 
union of points corresponding to the `$(l,l)$-isogenous abelian surfaces'. As there are
$[\Sp_4(\Z) : \Gamma_0^{(2)}(l)] = (l^4-1)/(l-1)$ isotropic subspaces of dimension
2 in $A[l]$ by~\cite[Lemma 6.1]{BrokerLauterPrePub}, there are exactly $(l^4-1)/(l-1)$ solutions to the system of
equations given by~$V'$. By construction, the triple
$$
(j_1(A/J), j_2(A/J), j_3(A/J))
$$
is one of the solutions. There are $l^3+l^2+l$ other solutions, and we will
see in Section~\ref{sec6} that for CM-computations it is relatively easy to determine 
which of the solutions come from the typenorm of an $\Order_{K_\Phi}$-ideal.

Unfortunately, the ideal $V(l)$ can only be computed for very small~$l$. Indeed,
the only case that has been done is $l=2$ (see~\cite[Sec.~10.4.2]{DupontThesis}) and it takes roughly 50 Megabytes
to store the 3 generators of $V$. By~\cite{BrokerLauterPrePub}, knowing the ideal $V(l)$ for some
prime $l$ implies that we have an equation for the {\it Humbert surface\/}
of discriminant~$l^2$. As computing Humbert surfaces is a traditionally hard
problem, we do not expect that much progress can be made in computing $V(l)$
for primes $l>2$. 

\subsection{Richelot isogeny}
Although one could use the ideal $V(2)$ from~\cite[Sec.~10.4.2]{DupontThesis}
to compute $(2,2)$-isogenies, there is a more efficient way. This alternative,
known as the {\it Richelot isogeny\/}, is classical and we recall it here
for convenience.
Let $K$ be a field of characteristic different from $2$, and let $C/K$ be a non-singular
genus $2$ curve. We can choose an equation $Y^2 = f(X)$ for $C$, with $f \in K[X]$
a monic polynomial of degree~$6$. Any factorization $f = ABC$ into three monic 
degree $2$ polynomials defines a genus $2$ curve~$C'$ given by
$$
\Delta Y^2 = [A,B][A,C][B,C],
$$
where $\Delta$ is the determinant of $A,B,C$ with respect to the basis $1,X,X^2$,
and $[A,B] = A'B-AB'$ with $A'$ the derivative of~$A$. This new curve is 
non-singular precisely when $\Delta$ is non-zero. 

One proves~\cite{BostMestre} that $C$ and~$C'$ are $(2,2)$-isogenous. It is not
hard to see that there are exactly $15 = (2^4-1)/(2-1)$ different curves~$C'$ that 
can be obtained this way. It follows that this construction
gives all $(2,2)$-isogenous Jacobians~$\Jac(C')$.

\section{Smaller functions} \label{sec4}
The Igusa functions introduced in Section~\ref{sec2} are `too large' to be practical
in our computation of the CM-action: currently we cannot compute an ideal describing
the variety $Y_0^{2}(l)$ for primes $l>2$. In this section we introduce smaller
functions $f_1,\ldots,f_4$ that are more convenient from a computational
perspective. For $N>1$, we define the {\it congruence subgroup\/} of level~$N$ as 
the kernel of the reduction map $\Sp_4(\Z) \rightarrow \Sp_4(\Z/N\Z)$, denoted by  $\Gamma(N)$. 

For $x,y \in \{0,1\}^2$, define the functions $\theta_{x,y}: \H_2 
\rightarrow \C$ by
\begin{equation*} 
\theta_{x,y}(\tau) = \sum_{n\in\Z^2} \exp \pi i\left( 
(n+\tfrac{x}{2})^T \tau (n+\tfrac{x}{2})+ (n+\tfrac{x}{2})^T y \right).
\eqno(4.1)
\end{equation*}
The functions $\theta_{x,y}$ are known as the {\it theta constants} and arise
naturally from the construction of theta functions \cite{IgusaThetaConsts1}. 
The 
equality $\theta_{x,y}(\tau) = (-1)^{x^Ty}\theta_{x,y}(\tau)$ shows that
only 10 of the 16 theta constants are non-zero.

The fourth powers of the functions $\theta_{x,y}$ are Siegel modular forms
of weight~2 for the congruence subgroup $\Gamma(2)\subset \Sp_4(\Z)$. The
Satake compactification $X(2)$ of the quotient $\Gamma(2) \backslash \H_2$
has a natural structure of a projective variety, and the fourth powers
$\theta_{x,y}^4$ define an embedding of $X(2)$ into projective space.

\begin{theorem} \label{sec4:thm1}
Let $M_2(\Gamma(2))$ denote the $\C$-vector space of
all Siegel modular forms of weight~2 for the congruence subgroup~$\Gamma(2)$. Then the
following holds: the space $M_2(\Gamma(2))$ is $5$-dimensional and is spanned
by the ten modular forms $\theta_{x,y}^4$. Furthermore, the map $X(2) 
\rightarrow \P^4 \subset \P^9$ defined by the functions $\theta_{x,y}^4$ is an embedding.
The image is the quartic threefold in $\P^4$ defined by
$$
u_2^2-4u_4=0
$$
with $u_k = {\displaystyle \sum_{x,y}\theta_{x,y}^{4k}}$.
\end{theorem}\noindent
{\Proof} See \cite[Thm.~5.2]{vdG82}. \endproof\par
\ \par\noindent

The Igusa functions $j_1,j_2,j_3$ can be readily expressed in terms of $\theta_{x,y}^4$, see e.g.~\cite[p.~848]{IgusaMPI}. Thus we have an inclusion
$$
\C(j_1,j_2,j_3) \subseteq \C(\theta_{x,y}^4/\theta_{x',y'}^4)
$$
where we include {\it all\/} quotients of theta
fourth powers. 
The functions $\theta^4_{x,y}/\theta^4_{x',y'}$ are rational Siegel modular
{\it functions\/} of level~$2$. Whereas 
$(j_1(\tau),j_2(\tau),j_3(\tau))$ depends only on the $\Sp_4(\Z)$-equivalence
class of $\tau\in\H_2$, a value $(\theta^4_{x,y}(\tau)/\theta^4_{x',y'}(\tau))_{x,x',y,y'}$
depends on the $\Gamma(2)$-equivalence class of~$\tau$. Since the affine
points of
$\Gamma(2) \backslash \H_2 \subset X(2)$ correspond to isomorphism
classes of pairs $(A,\{P_1,P_2,P_3,P_4\})$ consisting of a principally
polarized $2$-dimensional abelian variety $A$ together with a basis 
$\{P_1,P_2,P_3,P_4\}$ of the $2$-torsion, the functions $\theta^4_{x,y}/\theta^4_{x',y'}$ not only depend on the abelian variety in question but also on an ordering of its
$2$-torsion. For every isomorphism class $\Sp_4(\Z)\tau$ of abelian varieties,
there are $[\Sp_4(\Z):\Gamma(2)] = 720$ values for the tuple
$(\theta^4_{x,y}(\tau)/ \theta^4_{x',y'}(\tau))_{x,x',y,y'}$.
The functions $\theta^4_{x,y}/ \theta^4_{x',y'}$ are `smaller'
than the Igusa functions in the sense that their Fourier coefficients are
smaller. A natural idea is to get even smaller functions by considering the
quotients $\theta_{x,y}/\theta_{x',y'}$ themselves instead of their fourth
powers.

We define the four functions $f_1,f_2,f_3,f_4: \H_2 \rightarrow \C$ by
$$
f_1 = \theta_{(0,0),(0,0)} \quad f_2 = \theta_{(0,0),(1,1) }
 \quad f_3 = \theta_{(0,0),(1,0)} \qquad f_4 = \theta_{(0,0),(0,1)},
$$
with $\theta_{(x,y),(x',y')} = \theta_{x,y}/\theta_{x',y'}$.
We stress that the particular choice of the `theta constants' is rather
arbitrary, our only requirement is that we define $4$ different functions. The
three quotients $f_1/f_4,f_2/f_4,f_3/f_4$ are rational Siegel modular
functions. \par
\begin{theorem} \label{sec4:thm2}
If $\tau,\tau' \in \H_2$ satisfy
$(f_1(\tau),\ldots,f_4(\tau))=(f_1(\tau '),\ldots,f_4(\tau '))$, 
then we have 
$(j_1(\tau),j_2(\tau),j_3(\tau))=(j_1(\tau '),j_2(\tau '),j_3(\tau '))$.  
Furthermore, the quotients $f_1/f_4,f_2/f_4,f_3/f_4$ are invariant under the subgroup $\Gamma(8)$.
\end{theorem}
\noindent
{\bf Proof.}\ The vector space $M_2(\Gamma(2))$ is spanned by 
$\{f_1^4,\ldots f_4^4,g^4\}$ where $g = \theta_{(0,1),(0,0)}$. 
The relation in Theorem~\ref{sec4:thm1}, together with the 
five linear relations between the $\theta_{x,y}^4$ from Riemann's theta
formula \cite[p.~232]{IgusaThetaConsts1}, yield
that $g^4$ satisfies a degree $4$ polynomial $P$ over $L=\C(f_1,f_2,f_3,f_4)$.
The polynomial $P$ factors over $L$ as a product of the 2 irreducible 
quadratic polynomials
$$
P_-,P_+ = T^2 - (f_1^4-f_2^4+f_3^4-f_4^4)T + (f_1^2f_3^2\pm f_2^2f_4^2)^2.
$$
By looking at the Fourier expansions of $f_1,\ldots,f_4$ and $g$, we see
that $g^4$ is a root of~$P_-$. Hence, the extension
$L(g^4)/L$ is quadratic and generated by a root of $P_-$.

For each of the 2 choices of a root of $P_-$, the other $5$ fourth powers
of theta functions will be uniquely determined. Indeed, the fourth powers are
functions on the space $M_2(\Gamma(2))$ and this space is $5$-dimensional by
Theorem~\ref{sec4:thm1}. This means that we get a priori {\it two\/} Igusa triples
$(j_1,j_2,j_3)$ for every tuple $(f_1,f_2,f_3,f_4)$. However, a close
inspection of the formulas expressing the Igusa functions in terms of
theta fourth powers yields that these Igusa triples coincide. Hence, the
triple $(j_1,j_2,j_3)$ does not depend on the choice of a root of $P_-$. This 
proves the first statement in the theorem.

The second statement follows immediately from a result of Igusa. In
\cite[p.~242]{IgusaThetaConsts1}, he proves that the field $M$ generated by
{\it all\/} theta quotients is invariant under a group that contains 
$\Gamma(8)$. As the field $\C(f_1/f_4,f_2/f_4,f_3/f_4)$ is a subfield of $M$,
Theorem~\ref{sec4:thm2} follows.
\endproof\par
\ \par\noindent
As the functions $f_1/f_4,f_2/f_4,f_3/f_4$ are invariant under $\Gamma(8)$,
the moduli interpretation is that they depend on an abelian variety together
with a level $8$-structure. We let $\Stab(f)$ be the stabilizer of 
$f_1/f_4,f_2/f_4,f_3/f_4$ inside the symplectic group $\Sp_4(\Z)$. We have
inclusions
$$
\Gamma(8) \subset \Stab(f) \subset \Sp_4(\Z)
$$
and the quotient $Y(f) = \Stab(f) \backslash \H_2$ has a natural structure
of a quasi-projective variety by the Baily-Borel theorem~\cite{BB}. However, this
variety is not smooth. 

We let 
$$
\H_2^* = \{ \tau \in \H_2 \mid \tau \hbox{\ is not $\Sp_4(\Z)$-equivalent to
a diagonal matrix} \}
$$
be the subset of $\H_2$ of those $\tau$'s that do not correspond to a 
product of elliptic curves with the product polarization. The argument in 
\cite[Sec.\ 5]{Rung93} shows that $G = \Gamma(8) / \Stab(f)$ acts freely on 
$Y(8)$. By \cite[Ch.~2, Sec.\ 7]{Mum}, the quotient 
$Y(f)^* = \Stab(f) \backslash \H_2^*$ is a {\it smooth\/} variety.

\begin{lemma} The map $Y(f)^* \rightarrow Y(1)$ induced by the inclusion 
$\Stab(f) \rightarrow \Sp_4(\Z)$ has degree $23040=32 \cdot 720$.
\end{lemma} 

\Proof
The map factors as $Y(f)^* \rightarrow Y(2) \rightarrow Y(1)$ thus it suffices to determine the degrees of each part.
The degree of the map  $Y(f)^* \rightarrow Y(2)$ can be seen from the proof of Theorem \ref{sec4:thm2}: given a projective tuple $(f_4^4,f_2^4,f_3^4,f_4,g^4)$ representing a point $Q$ of $Y(2)$, over a splitting field there are $4^3=64$ projective tuples $(f_1,f_2,f_3,f_4)$ and exactly half of these satisfy $P_-=0$ and hence valid preimages of $Q$.  Thus $Y(f)^* \rightarrow Y(2)$ has degree $32$. The degree of $Y(2) \rightarrow Y(1)$ equals $[\Sp_4(\Z):\Gamma(2)]=720$.  This completes the proof.
\endproof\par
\ \par\noindent
The proof of Theorem~\ref{sec4:thm2} gives a means
of computing an Igusa triple $(j_1(\tau),j_2(\tau),j_3(\tau))$ from a tuple
$(f_1(\tau),\ldots,f_4(\tau))$. For convenience, we make this explicit in the
next subsection.

\subsection{Transformation formulae}
As in the proof of Theorem~\ref{sec4:thm2}, we let $g = \theta_{(0,1),(0,0)}$.
The function $g^4$ is a root of the quadratic polynomial~$P_-$. Given 
values $(f_1,f_2,f_3,f_4)$, we can pick any root of $P_-$ as a value for~$g^4$.
The functions $\{f_1^4,\ldots f_4^4,g^4\}$ form a basis of~$M_2(\Gamma(2))$.  
Now, define new functions $x_i$ by
\newcommand{\changethisbasis}[5]{#1 f_1^4 #2 f_2^4 #3 f_3^4 #4 f_4^4 #5 g^4}
\newcommand{\changethisbasisnog}[4]{#1 f_1^4 #2 f_2^4 #3 f_3^4 #4 f_4^4} 
\begin{align*}
x_1 &= \changethisbasis{-}{+2}{-}{+2}{+3} \,,\\
x_2 &= \changethisbasisnog{-}{+2}{-}{-} \,,\\
x_3 &= \changethisbasisnog{-}{-}{-}{+2} \,,\\
x_4 &= \changethisbasisnog{2}{-}{-}{-} \,,\\
x_5 &= \changethisbasisnog{-}{-}{+2}{-} \,,\\
x_6 &= \changethisbasis{2}{-}{+2}{-}{-3} \,.\\
\end{align*}
The $x_i$ are called level $2$ {\it Satake coordinate functions\/}. In terms of these functions we obtain a model for $X[2]$ embedded in $\P^5$ given by
\begin{eqnarray*}
s_1=0\,,\\
s_2^2-4s_4=0\,,
\end{eqnarray*}
where $s_k=\sum_{i=1}^6 x_i^k$ are the $k$-th power sums. 

The action of $\Sp_4(\Z)/\Gamma(2)$ on $x_i(\tau)$ is equivalent to that 
of $\Sym(\{x_1,\ldots,x_6\})$ permuting the coordinates.  Thus we can write 
level $1$ modular functions as symmetric functions of the $x_i$.  In 
particular, the {\it Igusa-Clebsch invariants\/} from Section~2.2
are given by 
\begin{align*}
I_2 &= \frac{5(48s_6 - 3s_2^3 - 8s_3^2)}{3(12s_5-5s_2s_3)} \,, \\
I_4 &= 3^{-1}s_2^2 \,, \\
I_6 &= 3^{-2}(3I_2I_4-2s_3) \,, \\
I_{10} &= 2^{-2}3^{-6}5^{-1}(12s_5-5s_2s_3) \,,
\end{align*}
from which we can compute absolute Igusa invariants $(j_1,j_2,j_3)$.

Conversely, if $(j_i(\tau))$ corresponds 
to the Jacobian of a curve, then we can compute a value 
for~$(f_1(\tau),\ldots,f_4(\tau))$ as follows. 
First we compute the Igusa-Clebsch invariants,  
then we apply the transformation
\medskip
\centerline{
\vbox{
\halign{$#$&$#$\hfill\cr
s_2&=3I_4\cr
\noalign{\smallskip}
s_3&=3/2(I_2I_4-3I_6)\cr
\noalign{\smallskip}
s_5&=5/12s_2s_3+3^5\cdot 5 I_{10}\cr
\noalign{\smallskip}
s_6&=27/16I_4^3+1/6s_3^2+3^6/2^2I_2I_{10},\cr}
}}
\medskip\noindent
after which we compute the level $2$ Satake coordinate functions as the roots $x_1,\ldots,x_6$ of the {\it Satake sextic polynomial\/}
$$
X^6-\frac{1}{2}s_2 X^4 - \frac{1}{3}s_3 X^3 + \frac{1}{16}s_2^2 X^2 + (\frac{1
}{6} s_2 s_3 - \frac{1}{5} s_5)X + (\frac{1}{96} s_2^3+\frac{1}{18}s_3^2-\frac{1}{6}s_6)
$$
with coefficients in $\Q(s_2,s_3,s_5,s_6)$. One choice for
$f_1^4,f_2^4,f_3^4,f_4^4$ is given by\par
\medskip
\centerline{
\vbox{
\halign{$#$&$#$\hfill\cr
f_1^4&=(-x_2 - x_3 - x_5)/3\cr
\noalign{\smallskip}
f_2^4&=(-x_3-x_4-x_5)/3\cr
\noalign{\smallskip}
f_3^4&=(-x_2-x_3-x_4)/3\cr
\noalign{\smallskip}
f_4^4&=(-x_2-x_4-x_5)/3.\cr}
}}
\medskip\noindent
Finally, we extract fourth roots to find values for $(f_1(\tau),\ldots,
f_4(\tau))$ satisfying $P_-=0$.
It is easy to find a solution to $P_-=0$: 
 if $(f_1,\ldots,
f_4)$ is not a solution, then $(\sqrt{-1}f_1,\ldots,
f_4)$ is a solution.

The coefficients of the Satake sextic polynomial are elements of $\ZZ[\frac{1}{2},\frac{1}{3},I_2,I_4,I_6,I_{10}]$.  In particular, this means that our transformation formulae are also valid over finite fields of characteristic greater than $3$.

\section{The CM-action and level structure} \label{sec5}
We let $\Stab(f)$ be the stabilizer of the three quotients $f_1/f_4, f_2/f_4,
f_3/f_4$ defined in Section~\ref{sec4}. By Theorem~\ref{sec4:thm2}, we have $\Gamma(8) \subseteq
\Stab(f)$. For a prime $l>2$, we now define
$$
Y(f;l)^* = (\Stab(f) \cap \Gamma_0^{(2)}(l)) \backslash \H_2^*
$$
which we view as an equality of Riemann surfaces. By the Baily-Borel theorem,
the space $Y(f;l)^*$ has a natural structure of a variety. Since we restricted
to $\H_2^*$, the variety is affine. Just
like in the case $l=1$ from Section~\ref{sec4}, $Y(f;l)^*$ is smooth.

The moduli interpretation of $Y(f;l)^*$ is the following.
Points are isomorphism classes of triples $(A,G,L)$, where $A$ is a 
principally polarized complex abelian surface, $G$ is a
$2$-dimensional isotropic subspace of $A[l]$ and $L$ is a level $8$-structure. The
notion of isomorphism is that $(A,G,L)$ and $(A',G',L')$ are isomorphic if
and only if there is an isomorphism of principally polarized abelian 
surfaces
$$
\varphi: A \rightarrow A'
$$
that satisfies $\varphi(G) = G'$ and $\varphi(L) = L'$.

\begin{lemma} The map $Y(f;l)^* \rightarrow Y(f)^*$ induced by the inclusion
map $(\Stab(f) \cap \Gamma_0^{(2)}(l)) \rightarrow \Stab(f)$ has degree
$(l^4-1)/(l-1)$ for primes $l>2$.
\end{lemma} 
\Proof
This is clear: the choice of a level $8$-structure $L$ is independent of the
choice of a subspace of the $l$-torsion for $l>2$.
\endproof\par
\ \par\noindent
Besides the map $Y(f;l)^* \rightarrow Y(f)^*$ from the lemma, we also have a map
$Y(f;l)^* \rightarrow Y(f)^*$ given by
$$
(A,G,L) \mapsto (A/G,L').
$$
Indeed, the isogeny $\varphi: 
A \rightarrow A/G$ induces an isomorphism
$$
A[8] \rightarrow (A/G)[8]
$$
and we have $L' = \varphi(L)$. As was explained in Section~3.2, this map also has
degree $(l^4-1)/(l-1)$. 
Putting all the varieties together, the picture is as follows.
$$
\xymatrix{& Y(f;l)^* \ar[dl]_s \ar[dr]^t \ar@/^.8pc/[ddl]_{} 
            \ar@/_.8pc/[ddr]^{}\\
Y(f)^* \ar[d]^{f} &&Y(f)^*\ar[d]_{f}\\
Y(1)&&Y(1)\ar@{.>}[d]\\
\A^3\ar@{->}[u]&&\A^3}
$$

The map $s$ sends $(A,G,L) \in Y(f;l)^*$ to $(A,L) \in Y(f)^*$ and $t$ is the
map induced by the isogeny $A \rightarrow A/G$. This diagram allows us
to find all the abelian surfaces that are $(l,l)$-isogenous to a given
surface~$A$, where we assume that $A$ is the Jacobian of a genus 2 curve. 
Indeed, we first map the Igusa invariants $(j_1(A),j_2(A),j_3(A))$
to a point in $Y(1)$, say given by the Igusa-Clebsch invariants. We then
{\it choose\/} $(A,L)$ on $Y(f)^*$ lying over this point. Although there are
$23040$ choices for $L$, it does not matter which one we choose. Above
$(A,L)$, there are $(l^4-1)/(l-1)$ points in $Y(f;l)^*$ and via the map
$t: Y(f;l)^* \rightarrow Y(f)^*$ we map all of these down to $Y(f)^*$. Forgetting the
level $8$-structure now yields $(l^4-1)/(l-1)$ points in $Y(1)$. If $A$ is
simple, i.e., not isogenous to a product of elliptic curves, then we can 
transform these into absolute Igusa invariants.

Assuming we can compute an ideal $V(f;l) 
\subset \Q[W_1,X_1,Y_1,Z_1,W_2,X_2,Y_2,Z_2]$ defining the quasi-projective 
variety $Y(f;l)^*$, we derive the following algorithm to compute all 
$(l,l)$-isogenous abelian surfaces.  

\begin{algorithm} \label{sec5:alg}
\ \par\noindent
{\bf Input.} A Jacobian $A/\C$ of a genus 2 curve given by its
Igusa invariants, and the ideal $V(f;l)$ defining $Y(f;l)^*$.\par\noindent
{\bf Output.} The Igusa invariants of all principally polarized abelian 
surfaces that are $(l,l)$-isogenous to $A$.
\begin{enumerate}
\item
Compute Igusa-Clebsch invariants $(I_2,I_4,I_6,I_{10})\in \C^4$ 
corresponding to~$A$.
\item 
Choose an element $(f_1,f_2,f_3,f_4) \in Y(f)^*$ that maps to 
$(I_2,I_4,I_6,I_{10})$
using the method described in Section~4.1.
\item 
Specialize the ideal $V(f;l)$ in $(W_1,X_1,Y_1,Z_1)=(f_1,f_2,f_3,f_4)$ and 
solve the remaining system of equations. 
\item
For each solution found in the previous step, compute the corresponding
point in $Y(1)$ using the method given in Section~4.1.
\end{enumerate}
\end{algorithm}

\subsection{Computing $V(f;l)$}\par
\ \par
In this subsection, we give an algorithm from~\cite{Gruenewald} to compute the ideal $V(f;l)$ needed in Algorithm~\ref{sec5:alg}. Our approach only terminates in a reasonable amount of time in the simplest case $l=3$. 

The expression for the theta constants in~(4.1) can be written in terms of the individual matrix entries, and with some minor modifications we can represent it as a power series with integer coefficients.
Write $\tau = \smat{\tau_1}{\tau_2}{\tau_2}{\tau_3} \in \H_2$, then
\begin{eqnarray*}
      \theta_{(a,b),(c,d)}(\tau)
&=& %
(-1)^{\frac{ac+bd}{2}}
\!\!\!\!\!\!\!
  \sum_{(x_1, x_2) \in \ZZ^2} (-1)^{x_1c+x_2d}  p^{(2x_1+a)^2}
  q^{(2x_1+a + 2x_2+b)^2} r^{(2x_2+b)^2} \in \Z[[p,q,r]]
\end{eqnarray*}
where $p=e^{2\pi i(\tau_1-\tau_2)/8}$, $q=e^{2\pi i\tau_2/8}$ and $r=e^{2\pi
  i(\tau_3-\tau_2)/8}$. 
We see that it is easy to compute Fourier expansions for the Siegel
modular forms~$f_i$.

One of the $(l,l)$-isogenous surfaces to $\C^2/(\Z^2+\Z^2\cdot\tau)$ is the
surface 
$$
\C^2/(\Z^2+\Z^2\cdot l\tau),
$$
and we want to find a relation between
the $f_i$'s and the functions $f_i(l\tau)$.
The expansion for
$f_i(l\tau)$ can be constructed easily from the Fourier expansion of 
$f_i(\tau)$ by replacing $p,q,r$ with $p^l,q^l,r^l$.

Starting with $n=2$, we compute all homogeneous monomials of degree $n$ in 
$f_i(\tau),f_i(l\tau)$ represented as truncated power series and then 
use exact linear algebra to find linear dependencies between them.
In this manner we obtain a basis for the degree $n$ homogeneous component of 
the relation ideal. We then check experimentally whether our list of relations 
generate $V(f;l)$ or not by computing the degree of the projection maps.
If one of the projection maps has degree larger than $l^3+l^2+l+1$, then more 
relations are required, in which case we increment $n$ by 1 and repeat the 
procedure. We stop once we have found sufficiently many relations to 
generate~$V(f;l)$.

Using this method we computed the ideal $V(f;3)$.
The $(3,3)$-isogeny relations in $V(f;3)$ are given by $85$ homogeneous 
polynomials of degree six.  
The whole ideal takes $35$ kilobytes to store in a text file, see\par
\centerline{
{\tt http://echidna.maths.usyd.edu.au/$\sim$davidg/ThesisData/theta33isogeny.txt.}
}
The individual relations are fairly small, having at most 40 terms. Furthermore, the coefficients are 
$7$-smooth and bounded by $200$ in absolute value, which makes them amenable 
for computations.

We point out however that we have not rigorously proven that the ideal $V(f;3)$ is correct.  
To do this we would need to show that our $85$ polynomials define relations between Siegel modular forms rather than just truncated Fourier expansions.
From the work of Poor and Yuen \cite{PoorYuen} there is a computable bound 
for which a truncated Fourier expansion uniquely determines the underlying Siegel modular form.  Thus with high enough precision our relations are able to be proven.
A Gr\"obner basis computation in Magma~\cite{Magma} informs us that the projection maps have the expected degree 40, hence we have obtained enough relations. 
Under the assumption that these relations hold, the ideal $V(f;3)$ is correct.

Our $(3,3)$-isogeny relations hold for all Jacobians of curves. We remark that
if we restrict ourselves to CM-abelian surfaces defined over unramified
extensions of~$\Z_3$, then there are smaller $(3,3)$-isogeny relations, 
see~\cite{3adic}. These smaller relations cannot be used however to improve the 
`CRT-algorithm' as in Section~6.3.

\section{The CM-action over finite fields} \label{sec6}
\subsection{Reduction theory}\par
\ \par
The theory developed in Sections \ref{sec3}--\ref{sec5} uses the {\it complex analytic\/} 
definition of abelian surfaces and the Riemann surfaces $Y_0^{(2)}(l)$ 
and $Y(f;l)^*$. We now explain why we can use the results in {\it positive 
characteristic\/} as well. Firstly, if we take a prime $p$ that splits
completely in $K$, then by \cite[Thm.\ 1, Thm.\ 2]{Gor97} the reduction modulo~$p$ of an abelian 
surface $A/H(K_\Phi)$ with endomorphism ring $\Order_K$ is {\it ordinary\/}. The 
reduced surface again has endomorphism ring~$\Order_K$.

Furthermore, one can naturally associate an algebraic stack 
$\sA_{\Gamma_0(p)}$ to $Y_0^{(2)}(l)$ and prove that the structural morphism
$\sA_{\Gamma_0(p)} \rightarrow \Spec(\Z)$ is smooth outside~$l$, see 
\cite[Cor.\ 6.1.1.]{CN}. In more down-to-earth computational terminology, 
this means the moduli interpretation of the ideal $V \subset 
\Q[X_1,\ldots,Z_2]$ remains valid when we reduce the elements of $V$ modulo a 
prime $p \not = l$. 

The reduction of $Y(f;l)^*$ is slightly more complicated. The map 
$Y(8l) \rightarrow Y(f;l)^*$ is finite \'etale by \cite[Thm.\ A.7.1.1.]{KM},
where we now view the affine varieties $Y(f;l)^*$ and $Y(8l)$ as schemes. 
It is well known that 
$Y(N)$ is smooth over $\Spec(\Z[1/N])$ for $N \geq 3$, so
in particular, the scheme $Y(f;l)^*$ is smooth over $\Spec(\Z[1/(2l)])$.
Again, this means that the
moduli interpretation for the ideal $V(f;l) \subset \Q[W_1,\ldots,Z_2]$
remains valid when we reduce the elements of $V(f;l)$ modulo a prime 
$p \nmid 2l$.  

We saw at the end of Section \ref{sec4} that our transformation formulae 
are valid modulo $p$ for primes $p>3$.
Putting this all together, we obtain the following result:

\begin{lemma} \label{sec6:lemma} Let $l$ be prime, and let $p \nmid 6l$ be a prime that splits
completely in a primitive CM-field $K$. Then, on input of the Igusa
invariants of a principally polarized abelian surface $A/\Fpbar$ with 
$\End(A) = \Order_K$ and the ideal $V(f;l) \subset \Fpbar[W_1,\ldots,Z_2]$, 
Algorithm~\ref{sec5:alg} computes the Igusa invariants of all $(l,l)$-isogenous abelian surfaces.
\end{lemma}

\subsection{Finding $(l,l)$-isogenous abelian surfaces}\par
\ \par
Fix a primitive quartic CM-field $K$, and let $p \nmid 6l$ be a prime 
that splits completely 
in the subfield 
$K_\Phi(j_1(A),j_2(A),j_3(A))$ of the Hilbert class field of $K_\Phi$.
By the choice of $p$, the Igusa invariants of an abelian surface $A/\Fpbar$ with $\End(A) = \Order_K$ are defined over the prime field $\Fp$.
Moreover, $p$ splits in $K_\Phi$ and as it splits
in its normal closure~$L$ it will split completely in~$K$, 
hence Lemma~\ref{sec6:lemma} applies.

If we apply Algorithm~\ref{sec5:alg} to the point $(j_1(A),j_2(A),j_3(A))$ and the ideal $V(f;l)$,
then we get $(l^4-1)/(l-1)$ triples of Igusa invariants.  All these triples
are Igusa invariants of principally polarized abelian surfaces with
endomorphism {\it algebra\/}~$K$. Some of these triples are defined over the
prime field $\Fp$ and some are not. However, since $p$ splits completely in
the field of moduli $K_\Phi(j_1(A),j_2(A),j_3(A))$, the Igusa invariants of the surfaces 
that have endomorphism ring $\Order_K$ {\it are defined over the field\/} $\Fp$.

\begin{algorithm} \label{sec6:cm}
\ \par\noindent
{\bf Input.} The Igusa invariants of a simple principally polarized abelian
surface $A/\Fp$ with $\End(A) = \Order_K$, and the ideal $V(f;l) \subset
\Fp[W_1,\ldots,Z_2]$. Here, $l$ is a prime such that there exists a
prime ideal in $K_\Phi$ of norm $l$. Furthermore, we assume that $p \nmid 6l$.
\par\noindent
{\bf Output.} The Igusa invariants of all principally polarized abelian 
surfaces $A'/\Fp$ with $\End(A') = \Order_K$ that are $(l,l)$-isogenous to $A$.
\begin{enumerate}
\item
Apply Algorithm~\ref{sec5:alg} to $A$ and $V(f;l)$. Let $S$ be the set of all Igusa
invariants that are defined over $\Fp$.
\item 
For each $(j_1(A'),j_2(A'),j_3(A')) \in S$, construct a genus 2 curve $C$
having these invariants using Mestre's algorithm~(\cite{Mestre},~\cite{CardonaQuer}).
\item 
Apply the Freeman-Lauter algorithm~\cite{FreemanLauter} to test whether $\Jac(C)$ has
endomorphism ring $\Order_K$. Return the Igusa invariants of all the curves that
pass this test.
\end{enumerate}
\end{algorithm}
\noindent
We can predict beforehand {\it how many\/} triples will be returned by 
Algorithm~\ref{sec6:cm}. We compute the prime factorization
$$
(l) = \gp_1^{e_1} \ldots \gp_k^{e_k} 
$$
of $(l)$ in $K_\Phi$. Say that we have $n \leq 4$ prime ideals $\gp_1,\ldots, \gp_n$ 
of norm~$l$ in this factorization, disregarding multiplicity. For each of 
these ideals $\gp_i$ we compute the typenorm map $m(\gp_i) \in \gC(K)$. 
The size of 
$$
\{ m(\gp_1),\ldots,m(\gp_n)\} \subset \gC(K).
$$
equals the number of triples computed by Algorithm~\ref{sec6:cm}. 

\begin{remark}
Step~1 of the algorithm requires working in an extension of~$\Fp$. The 
degree of this extension depends on the splitting behaviour of~$2$ 
in~$\Order_K$. An upper bound is given by $4 [\Fp(A[2]) : \Fp] \leq 24$, 
where $\Fp(A[2])$ denotes the field obtained by adjoining the coordinates of 
all $2$-torsion points of~$A$.
\end{remark}

\subsection{Igusa class polynomials}\label{sec6.3} \par
\ \par
The CRT-algorithm~\cite{EisentragerLauter} to compute the Igusa class 
polynomials $P_K, Q_K, R_K \in \Q[X]$ of a primitive quartic CM-field~$K$ 
computes the reductions of these $3$ polynomials modulo primes~$p$ which split 
completely in the Hilbert class field of $K_\Phi$. The method suggested 
in~\cite{EisentragerLauter} loops over all $p^3$ possible Igusa invariants and 
runs an endomorphism ring test for each triple  $(j_1(A'),j_2(A'),j_3(A'))$, 
to see if $A'$ has endomorphism ring~$\Order_K$. 

We propose two main modifications to this algorithm. Firstly, we only
demand that the primes~$p$ split completely in the subfield 
$K_\Phi(j_1(A),j_2(A),j_3(A))$ of the Hilbert class field of $K_\Phi$ that
we obtain by adjoining the Igusa invariants of an abelian surface~$A$ that
has CM by~$\Order_K$. To find such primes, we simply loop over $p=5,7,11\ldots$,
and for the primes $(p) = \gP_1\gP_2\gP_3\gP_4 \subset \Order_{K_\Phi}$ that 
split completely in~$K_\Phi$ we test if
$$
m(\gP_1) = (\mu) \subset \Order_K \qquad \hbox{\ and \ } N(\gP_1) = \mu 
\overline \mu
$$
holds, with $\overline\mu$ the complex conjugate of~$\mu$. 
By~\cite[Sec.~15.3, Thm.~1]{Shim} a prime $p$ satisfying these conditions
splits completely in~$K_\Phi(j_1(A),j_2(A),j_3(A))$. It includes the primes
that split completely in the Hilbert class field of~$K_\Phi$.

Our second modification is a big improvement to compute the Igusa class
polynomial modulo~$p$. Instead of looping over all $O(p^3)$ curves, we
exploit the Galois action in a similar vein as in~\cite{BBEL}. Below we
give the complete algorithm to compute $P_K,R_K,Q_K$ modulo a prime
$p$ meeting our conditions.\par
\ \par\noindent
{\bf Step 1.} Compute the class group 
$$
\Cl(\Order_{K_\Phi}) = \langle \gp_1,\ldots,\gp_k \rangle \eqno(6.1)
$$
of the reflex field, where we take degree 1 prime ideals~$\gp_i$. 
For each of 
the ideals $\gp_i$ of {\it odd\/} norm $N_{K_\Phi/\Q}(\gp_i) = l_i$, 
compute the ideal $V(f;l_i)$ describing the Siegel modular variety $Y(f;l_i)^*$.\par
\ \par\noindent
{\bf Step 2.} Apply the following algorithm to find an abelian surface $A/\Fp$
that has endomorphism ring isomorphic to~$\Order_K$. We factor $(p) \subseteq
\Order_K$ into primes $\gP_1,\overline \gP_1,\gP_2,\overline \gP_2$ and
compute a generator $\pi$ for the principal ideal~$\gP_1 \gP_2$. We compute
the minimal polynomial~$f_\pi$ of $\pi$ over~$\Q$. We try {\it random\/} curves $C/\Fp$ until we 
find a curve with
$$
\#C(\Fp) \in \{p+1\pm \Tr_{K/\Q}(\pi)\} \quad \hbox{\ and \ } \#\Jac(C) 
\in \{ f_\pi(1), f_\pi(-1)\} \,. \eqno(6.2)
$$
By construction, such a curve~$C$ has endomorphism algebra $K$.  
We test whether $\Jac(C)$ has endomorphism ring~$\Order_K$ using
the algorithm from~\cite{FreemanLauter}. If it does, continue with Step~3,
otherwise try more random curves~$C$ until we find one for which its 
Jacobian has endomorphism ring~$\Order_K$.\par
\ \par\noindent
{\bf Step 3.} Let $A/\Fp$ be the surface found in Step~2.
The group $G = m(\Cl(\Order_{K_\Phi}))$ acts in a 
natural way on $A$ and we compute the set
$$
G \cdot (j_1(A),j_2(A),j_3(A)) \subseteq S(K)
$$
as follows. For $x = m(I)$ we write $I = \prod_i \gp_i^{a_i}$. The action 
of $\gp_1$ is computed using 
Algorithm~\ref{sec6:cm} in case the norm of $\gp_1$ is odd and by applying
a Richelot isogeny (see Section~3.3) if $\gp_1$ has norm~2. By successively 
applying the action of $\gp_1$ we compute the action of $\gp_1^{a_1}$. We
then continue with the action of~$\gp_2$, etc. This allows us to compute
the action of~$x$ on the surface~$A$, and doing this for all $x$ we 
compute the set $G \cdot (j_1(A),j_2(A),j_3(A))$. This part of the 
algorithm is analogous to~\cite{BBEL}.\par
\ \par\noindent
{\bf Step 4.} 
In contrast to genus $1$ and the algorithm in~\cite{BBEL}, it is unlikely 
that \textit{all} surfaces with endomorphism ring $\Order_K$ are found. This 
is partly because we only find surfaces having the \textit{same} CM-type as the 
initial surface $A$, so in the dihedral case we are missing surfaces with the 
second CM-type. Even in the cyclic case where there are $|\gC(K)|$ isomorphism 
classes, it is possible that the map
$$
m: \Cl(\Order_{K_\Phi}) \rightarrow \gC(K)
$$
is not surjective, meaning that we do not find all surfaces of a given CM-type.
The solution is simple: compute the cardinality of
$S(K)$ using Corollary~\ref{sec2:cardinality},
 and if the number of surfaces found is less than $|S(K)|$, go back to
Step~2.\par
\ \par\noindent
{\bf Step 5.}
Once we have found all surfaces with endomorphism ring $\Order_K$, 
expand
$$
P_K \bmod p= \prod_{\{A\ p.p.a.s. \mid \End(A) = \Order_K\}/\cong} (X - j_1(A)) \in \Fp[X]
$$
and likewise for $Q_K$ and $R_K$. The main difference with the method
from~\cite{EisentragerLauter} is that we do {\it not\/} find all roots 
of~$P_K$ by a random search: we exploit the Galois action. 

\subsection{Run time analysis}
We proceed with the run time analysis of the algorithm to compute the
Igusa class polynomials using the `modified CRT-approach' from Section~\ref{sec6.3}.
The input of the algorithm is a degree four CM-field~$K$. The discriminant~$D$
of~$K$ can be written as $D_1 D_0^2$, with $D_0$ the discriminant of the
real quadratic subfield~$K^+$ of~$K$. We will give the run time in terms
of~$D_1$ and~$D_0$.

First we analyze the size of the primes~$p$ used in the algorithm. The 
primes we use split completely in a subfield $S = K_\Phi(j_i(A))$ of the 
Hilbert class field of 
the reflex field~$K_\Phi$ of~$K$. If GRH holds true, then there exists~\cite{LO}
an effectively computable constant $c>0$, independent of~$K$, such that the
smallest such prime $p$ satisfies
$$
p \leq c\cdot (\log |\disc(S/\Q)|)^2,
$$
where $\disc(S/\Q)$ denotes the discriminant of the extension $S/\Q$. Since
$S$ is a totally unramified extension of~$K_\Phi$, we have
$$
\disc(S/\Q)^{1/[S:\Q]} = \disc(K_\Phi/\Q)^{1/[K_\Phi:\Q]}
$$
and we derive $\disc(S/\Q) = \disc(K_\Phi/\Q)^{[S:K_\Phi]}$. Theorem~\ref{sec2:thm3} 
yields the bound $[S:K_\Phi] \leq 4 h^-(K)$, where $h^-(K) = |\Cl(\Order_K)|/|\Cl(\Order_{K^+})|$ denotes the
{\it relative\/} class number of~$K$. Using the bound (see~\cite{Louboutin})
$$
h^-(K) = \Otilde (\sqrt{D_1 D_0}) \,, \eqno(6.3)
$$
we derive that the smallest prime~$p$ is of size~$\Otilde(D_1 D_0)$. Here,
the $\Otilde$-notation indicates that factors that are of logarithmic order
in the main term have been disregarded.

The Igusa class polynomials have rational coefficients, and at the moment
the best known bound for the logarithmic height of the denominator of a 
coefficient is $\Otilde(D_1^{3/2} D_0^{5/2})$. This bound is proven 
in~\cite[Sec.\ 2.9]{Streng} and is based on the denominator bounds 
in~\cite{GorenLauter}. A careful
analysis~\cite[Sec.\ 2.11]{Streng} yields that each coefficient of $P_K,
R_K, Q_K$ has logarithmic
height~$\Otilde(D_1^{3/2} D_0^{5/2})$ as well. A standard argument
as in~\cite[Lemma 5.3]{BBEL} shows that the $\Otilde(D_1^{3/2} 
D_0^{5/2})$ primes that we need can be taken to be of size 
$\Otilde(D_1^{2} D_0^{3})$ if GRH holds true. We find these primes in 
time $\Otilde(D_1^{2} D_0^{3})$. We remark that better bounds
on the denominators of the coefficients translate into better bounds on 
the size of the primes we need. 

If GRH holds true, then the ideals $\gp_i$ in Step~1 can be chosen to
have norm 
at most $12 (\log {D_1 D_0^2})^2$ by~\cite{Bach}. As the method from 
Section~\ref{sec5} to compute the ideal $V(f;N(\gp_i))$ is heuristic, we will rely
on the following heuristic for our analysis.

\begin{heuristic} \label{sec6:heu1}
Given a prime $l>2$, we can compute generators for the
ideal $V(f;l)$ in time polynomial in~$l$.
\end{heuristic}
\noindent
We remark that, at the moment, our implementation of computing~$V(f;l)$ only
terminates in a reasonable amount of 
time for~$l=3$. However, {\it in theory\/} we only spend heuristic time 
$(\log {D_1 D_0^2})^n$ in Step~1 for some~$n\geq 2$ that is independent of
$D_1,D_0$. This is negligible compared to other parts of the algorithm.

We continue with the analysis of computing $P_K \bmod p$. As we think that
the bound $p = \Otilde(D_1^{2} D_0^{3})$ is too pessimistic, we will
do the analysis in terms of both $p$ and $D_1,D_0$. First we analyse the
time spent on the random searches to find abelian surfaces with 
endomorphism ring~$\Order_K$. Every time we leave Step~2, we compute a 
factor $F \mid P_K \bmod p$ of the (first) Igusa class polynomial. 
Let $k \leq 2 [\gC(K) : m(\Cl(\Order_{K_\Phi}))]$
be the number of factors~$F$ we need to compute. The first time we invoke
Step~2, we will with probability~1 compute a new factor~$F_1$ of~$P_K$. The
second time we call Step~2 we need to ensure that we compute a {\it
different\/} factor~$F_2 \mid P_K$. Hence, we expect that we need to call
Step~2 $k/(k-1)$ times to compute~$F_2$. We see that we expect that we have
to do Step~2 
$$
k(1+1/2+\ldots+1/k) = \Otilde(k)
$$
times to compute all factors~$F_1,\ldots,F_k$.

\begin{lemma}\label{sec6:lemgenus} We have $[\gC(K) : m(\Cl(\Order_{K_\Phi}))] \leq 2 \cdot 
2^{6\omega(D)}$ for any primitive quartic CM-field~$K$, where $\omega(D)$ denotes
the number of prime divisors of~$D$.
\end{lemma}
\noindent
{\bf Proof.} We will bound the index of the image of the map 
$\tilde m: \Cl(\Order_{K_\Phi}) \rightarrow \Cl(\Order_K)$ inducing~$m$. By 
Theorem~\ref{sec2:thm3}, this index differs by at most a factor
$$ 
|(\Order_{K^+}^*)^+/N_{K/K^+}(\Order_K^*)| \leq 2
$$
from $[\gC(K) : m(\Cl(\Order_{K_\Phi}))]$.

If $K/\Q$ is dihedral with normal closure~$L$, then the image of the 
norm map $N_{L/K}: \Cl(\Order_L) \rightarrow \Cl(\Order_K)$ has index
at most 2 by class field theory. In 
the cyclic
case, it is not hard to check that $\Im(\tilde m)$ contains the squares. It
suffices to bound the 2-torsion $\Cl(\Order_{K})$ in this case. The 2-rank
of $\Cl(\Order_K)$ is determined by {\it genus theory\/}. Using a combination
of group cohomology and Nakayama's lemma, one can show~\cite{MR} that the
2-rank is at most $6t$, with $t$ the number of primes that ramify in the
cyclic CM-extension~$K/\Q$. The lemma follows.
\endproof\par
\ \par\noindent
We remark 
that outside a zero-density subset of very smooth integers, we have 
$\omega(n) < 2 \log\log n$ and the factor~$\Otilde(2^{6\omega(D)}) 
= 2^{6\omega(D)} \Otilde(\log(D))$ can then be absorbed into 
the $\Otilde$-notation.

The probability that one of the random searches performed in this step will yield 
an abelian surface $\Jac(C)$ with endomorphism ring $\Order_K$ is 
bounded from below by
$$
h^-(K)/p^3 = \tilde\Omega (\sqrt{D_1 D_0}/p^3) 
$$
where we have used the effective lower bound~$h^-(K) = \tilde\Omega(\sqrt{D_1 D_0})$
proved in~\cite{Louboutin}. We therefore expect that we have to compute the number
of points on $C$ and on $\Jac(C)$ for 
$$
\Otilde(p^3 / \sqrt{D_1 D_0})
$$
curves~$C/\Fp$. Since point counting on genus 2 curves is polynomial time by~\cite{Pila}, 
this takes time $\Otilde(p^3/\sqrt{D_1 D_0})$. 

For all the curves~$C/\Fp$ that satisfy equation~(6.2), we
have to check whether we have $\End(\Jac(C)) \cong \Order_K$ or not. The probability
that $\End(\Jac(C))$ is isomorphic to $\Order_K$ is bounded from below by
$$
\frac{h^-(K)}{\sum_\Order h(\Order)}, \eqno(6.4)
$$
where the sum ranges over all orders $\Order \subseteq \Order_K$ that contain
$\Z[\pi,\overline\pi]$. Assuming mild ramification conditions on the prime 2,
there are only $O(\log n)$ orders $\Order \subseteq \Order_K$ of index~$n$,
see~\cite[Cor.~1]{JN}. We assume the following heuristic.

\begin{heuristic} \label{sec6:heu3}
For any quartic CM-field~$K$, there are $O(\log n)$ orders $\Order \subseteq
\Order_K$ of index~$n$.
\end{heuristic}
\noindent
{\bf Justification of heuristic.} As indicated in~\cite{JN}, the splitting
condition on $2$ is purely technical and should not affect the result. 
\endproof\par
\ \par\noindent
We can bound the class number $h(\Order)$ by $2 
[\Order_K : \Z[\pi,\overline\pi]]h(\Order_K)$ by~\cite[Thm.~6.7]{psh}. 
It follows that we can bound the probability in~(6.4) 
by 
$$
\Omega\left(\frac{1}{[\Order_K : \Z[\pi,\overline\pi]]^{1+\varepsilon}
h(\Order_{K^+})}\right),
$$
where we have used the bound $n^\varepsilon$ for the number of divisors of~$n$.
Using the index bound $[\Order_K: \Z[\pi,\overline\pi]] \leq \frac{16p^2}{D_0\sqrt{D_1}}$ from~\cite[Prop.~6.1]{FreemanLauter}, we expect that
we have to do 
$$
O\left(\frac{ 2^{6\omega(D)}p^{2+2\varepsilon}\sqrt{D_0} }{ (\sqrt{D_1} D_0)^{1+\varepsilon} } \right)
$$
endomorphism ring computations. 

At the moment, the only known algorithm~\cite{FreemanLauter}
to test whether $\End(\Jac(C)) \cong \Order_K$ holds has a run time~$\Otilde(p^{18})$,
and one application of this algorithm dominates the computation of~$P_K \in 
\Q[X]$. To make the run time analysis of our algorithm easier once a better
algorithm to compute $\End(\Jac(C))$ has been found, we will use the bound $O(X)$
for the run time to compute~$\End(\Jac(C))$. In total, we see that we spend
$$
\Otilde(2^{6\omega(D)}(p^3 / \sqrt{D_1 D_0} + (p^{2+2\varepsilon} \sqrt{D_0}/ (\sqrt{D_1} D_0)^{1+\varepsilon})X))
$$
time in all calls of~Step~2.

The action of $\gp_i$ on $A/\Fp$ in Step~3 is computed in polynomial time
in the norm $l_i$ of~$\gp_i$. As $l_i$ is, under GRH, of polynomial size in 
$\log(D_1 D_0^2)$, we spend time 
$$
\Otilde(\sqrt{D_1 D_0})
$$
for every time we call Step~3. We call Step~3 as often as Step~2, so in total we
spend time $\Otilde(2^{6\omega(D)}\sqrt{D_1 D_0})$ in Step~3.

The time spent in Step~4 is negligible, and the time spent in Step~5 is
$\Otilde(\sqrt{D_1 D_0})$. Combining all five steps, we see that we compute
$P_K \bmod p$ in time
$$
\Otilde(2^{6\omega(D)}(p^3 / \sqrt{D_1 D_0} + (p^{2+2\varepsilon} \sqrt{D_0}/ (\sqrt{D_1} D_0)^{1+\varepsilon})X+ \sqrt{D_1 D_0})). \eqno(6.5)
$$

\begin{theorem} If GRH and Heuristic assumptions \ref{sec6:heu1}, 
\ref{sec6:heu3} hold true, then
we can compute the polynomials $P_K, Q_K, R_K$ in probabilistic time
$$
\Otilde(2^{6\omega(D)}(D_1^{7}D_0^{11}+X D_1^{5+\varepsilon} D_0^{8+2\varepsilon})).
$$
Here, $X$ denotes the run time of an algorithm that, given $A / \Fp$,
decides whether $\End(A)$ is isomorphic to $\Order_K$ or not.
\end{theorem}
\noindent
{\bf Proof.} Substitute $p = \Otilde(D_1^{2} D_0^{3})$ in 
equation~(6.5) to get the time per prime. The result follows from the fact
that we need to compute $P_K,R_K,Q_K$ modulo~$p$ 
for $\Otilde(D_1^{3/2} D_0^{5/2})$ primes.\endproof\par
\ \par\noindent
We conclude this section by some remarks on the run time of our algorithm. At
the moment, the main bottleneck in the run time is checking whether 
$\End(A) \cong \Order_K$
holds or not. In Section~\ref{sec8} we show that a `straightforward generalization'
of Kohel's algorithm~\cite{Kohel} is impossible and that a new approach is
needed.

From a practical point of view, we are limited by the fact that we can 
only compute the ideal~$V(f;3)$ in a reasonable amount of time. By only
using the primes lying over $2$ and $3$, we only use a subgroup of the 
group~$\gC(K)$ giving the Galois action. 

Even when these two problems are solved, there is a bottleneck that does not
occur in the genus $1$ algorithm from~\cite{BBEL}. The random searches take
time $\Otilde(p^3 / \sqrt{D_1 D_0})$ and even for the smallest prime~$p$
this is of size
$
\Otilde(D_1^{5/2} D_0^{5/2}).
$

Doing only the random searches for this prime already takes more time than
it takes to write down the output $P_K,Q_K,R_K \in \Q[X]$. Hence, our
algorithm is at the moment {\it not\/} quasi-linear in the size of the
output. 

As noted in~\cite[Sec.\ 6]{GruenewaldArticle}, 
we can speed up this step of the
algorithm by first computing a model
for the Humbert surface describing all principally polarized abelian 
surfaces that have real multiplication by the quadratic subfield~$K^+$
of~$K$. We then perform our random search on this two-dimensional subspace
of the three-dimensional moduli space. The time for the random searches would,
for the smallest prime~$p$, drop to
$$
\Otilde(D_1^{3/2} D_0^{3/2}).
$$
Although this is less than the size of the output, our algorithm is not
quasi-linear once all primes~$p$ are taken into account.

To get a quasi-linear algorithm, we think one should do the random searches
on a one-dimensional subspace of the moduli space.
This approach is an object of further study.

\section{Examples and applications} \label{sec7}

In this section we illustrate our algorithm by computing the Igusa class
polynomials modulo primes~$p$ for various CM-fields. We point out the
differences with the analogous genus 1 computations. 

\begin{example} \label{sec7:ex1}
In the first example we let $K= \Q[X]/(X^4 + 185X^2 + 8325)$ be a 
{\it cyclic\/} CM-field of degree~4. All CM-types are equivalent in this
case, and the reflex field of $K$ is $K$ itself. 
The discriminant of $K$ equals
$5^2 \cdot 37^3$, and the real quadratic subfield of $K$ is $K^+ = 
\Q(\sqrt{37})$. An easy computation shows that the narrow class group of
$K^+$ is trivial. In particular, all ideal classes of $K$ are principally
polarizable, and we have
$$
\mathfrak{C}(K) \cong \Cl(\Order_K).
$$
We compute $\Cl(\Order_K) = \Z/10\Z = \langle \gp_3 \rangle$, where $\gp_3$ is
a prime lying over $3$.
The prime ideal $\gp_3$ has norm $3$, and its typenorm
$N_\Phi(\gp_3)$ generates a subgroup of 
order~5 in $\Cl(\Order_K)$.

The smallest prime that splits in the Hilbert class field of $K$ is $p=271$. 
We illustrate our algorithm by computing the Igusa class polynomials for $K$
modulo this prime. First we do a `random search' to find a principally
polarized abelian surface over $\Fp$ with endomorphism ring $\Order_K$ in the
following way. We factor $(p) \subset \Order_K$ into primes $\gP_1,\overline \gP_1,
\gP_2,\overline\gP_2$ and compute a generator $\pi$ of the principal 
$\Order_K$-ideal $\gP_1 \gP_2$. The element $\pi$ has minimal 
polynomial
$$
f = X^4 + 9X^3 + 331X^2 + 2439X + 73441 \in \Z[X].
$$
If the Jacobian $\Jac(C)$ of a hyperelliptic curve $C$ 
has endomorphism ring $\Order_K$, then the Frobenius morphism of $\Jac(C)$
is a root of either $f(X)$ or $f(-X)$. With the factorization
$$
f = (X-\tau_1)(X-\tau_2)(X-\tau_3)(X-\tau_4) \in K[X],
$$
a {\it necessary\/} condition for $\Jac(C)$ to have endomorphism ring $\Order_K$
is 
$$
\#C(\Fp) \in \{ p+1 \pm (\tau_1+\tau_2+\tau_3+\tau_4)\} = \{261,283\}
$$
and
$$
\#\Jac(C)(\Fp) \in \{ f(1),f(-1)\} = \{71325, 76221\}.
$$
We try random values $(j_1,j_2,j_3)\in\Fp^3$ and write down a hyperelliptic
curve $C$ with those Igusa invariants using Mestre's algorithm~(\cite{Mestre},~\cite{CardonaQuer}). If
$C$ satisfies the 2 conditions above, then we check whether $\Jac(C)$ has
endomorphism ring $\Order_K$ using the algorithm in~\cite{FreemanLauter}. If it
passes this test, we are done. Otherwise, we select a new random value
$(j_1,j_2,j_3)$.

We find that $w_0 = (133,141,89)$ is a set of invariants for a surface
$A/\Fp$ with endomorphism ring $\Order_K$. We apply Algorithm~\ref{sec6:cm} to 
$w_0$. The Igusa-Clebsch invariants corresponding to $w_0$ are 
$[ 133, 54, 82, 56 ]$. With the notation from Section~\ref{sec4}, we have 
$s_2 = 162, s_3 = 106, s_5=128, s_6 = 30$. The Satake sextic polynomial
$$
\mathcal{S} = X^6 + 190X^4 + 55X^3 + 82X^2 + 18X + 63 \in \Fp[X]
$$
factors over $\F_{p^5}$ and we write $\F_{p^5} = \Fp(\alpha)$ where 
$\alpha$ satisfies $\alpha^5 + 2\alpha + 265 = 0$. We express
the 6 roots of $\mathcal{S}$ in terms of $\alpha$ and pick
\begin{align*}
f_1^4&= 147\alpha^4 + 147\alpha^3 + 259\alpha^2 + 34\alpha + 110\\
f_2^4&= 176\alpha^4 + 211\alpha^3 + 14\alpha^2 + 134\alpha + 190\\
f_3^4&= 163\alpha^4 + 93\alpha^3 + 134\alpha^2 + 196\alpha + 115\\
f_4^4&= 226\alpha^4 + 261\alpha^3 + 99\alpha^2 + 9\alpha + 27
\end{align*}
as values for the fourth powers of our Siegel modular functions. 
The fourth roots of $(f_1^4,f_2^4,f_3^4,f_4^4)$ are all defined over 
$\F_{p^{10}}$, but the proof of Theorem~4.2 shows that not every choice 
corresponds to the
Igusa invariants of $A$. We pick fourth 
roots $(r_1,r_2,r_3,r_4)$ such that the polynomial $P_-$ from 
Section~\ref{sec4} vanishes when evaluated at 
$(T,f_1,f_2,f_3,f_4) = (\theta_{(0,1),(0,0)}^4,r_1,r_2,r_3,r_4)$. Here, $\theta_{(0,1),(0,0)}^4$ is computed 
from the Igusa-Clebsch invariants. For an arbitrary choice of fourth roots 
for $r_1,r_2,r_3$ there are two solutions $\pm r_4$ to $P_-=0$. 
Indeed, if we take $\F_{p^{10}}=\Fp(\beta)$ with
$
\beta^{10} + \beta^6 + 133\beta^5 + 10\beta^4 + 256\beta^3 + 74\beta^2 + 126\beta + 6 =0
$
then the tuple $(r_1,r_2,r_3,r_4)$ given by
\footnotesize
\begin{align*}
r_1 &= 179\beta^9 + 69\beta^8 + 203\beta^7 + 150\beta^6 + 29\beta^5 + 258\beta^4 + 183\beta^3 + 240\beta^2 + 255\beta + 226 ,\\
r_2 &= 142\beta^9 + 105\beta^8 + 227\beta^7 + 244\beta^6 + 72\beta^5 + 155\beta^4 + 2\beta^3 + 129\beta^2 + 137\beta + 23 ,\\
r_3 &= 63\beta^9 + 112\beta^8 + 132\beta^7 + 244\beta^6 + 94\beta^5 + 40\beta^4 + 191\beta^3 + 263\beta^2 + 85\beta + 70 ,\\
r_4 &= 190\beta^9 + 41\beta^8 + 62\beta^7 + 170\beta^6 + 151\beta^5 + 240\beta^4 + 270\beta^3 + 56\beta^2 + 16\beta + 257
\end{align*}
\normalsize
is a set of invariants for $A$ together with some level $8$-structure.

Next we specialize our ideal $V(f;3)$ at $(W_1,X_1,Y_1,Z_1) = 
(r_1,r_2,r_3,r_4)$ and we solve the remaining system of 85 equations in
4 unknowns. Let $(r_1',r_2',r_3',r_4')$ be the solution where 
\small
\begin{align*}
r_1' &= 184\beta^9 + 48\beta^8 + 99\beta^7 + 83\beta^6 + 20\beta^5 + 232\beta^4 + 16\beta^3 + 223\beta^2 + 85\beta + 108 .
\end{align*}
\normalsize
The quadruple
$(r_1',r_2',r_3',r_4')$ are invariants of an abelian surface $A'$ together
with level $8$-structure that is $(3,3)$-isogenous to $A$. To map this
quadruple to the Igusa invariants of $A'$ we compute a root of the
quadratic polynomial 
$$
P_-(T,r_1',r_2',r_3',r_4').
$$
This root is a value for $\theta_{(0,1),(0,0)}^4$. Since we now
know {\it all\/} theta fourth powers, we can apply the formulas relating
theta functions and Igusa functions in Section~4.1 to find the Igusa triple
$( 238, 10, 158 )$.

In total, we find 16 Igusa triples defined over $\Fp$. All these triples
are Igusa invariants of surfaces that have endomorphism {\it algebra\/}~$K$.
To check which ones have endomorphism {\it ring\/}~$\Order_K$, we apply the
algorithm from~\cite{FreemanLauter}. We find that only the four triples
$$
(253,138,96),\ (257,248,58),\ (238,10,158),\ (140,159,219)
$$
are invariants of surfaces with endomorphism ring $\Order_K$. The fact that we
find 4 new sets of invariants should come as no surprise. Indeed, there are 4
ideals of norm 3 lying over $3$ in $\Order_K$ and each ideal gives us an
isogenous surface. 

As the typenorm map $m: \Cl(\Order_K) \rightarrow \mathfrak{C}(K)$ is not
surjective, we are forced to do a {\it second\/} random search to find a `new'
abelian surface with endomorphism ring $\Order_K$. We apply our isogeny 
algorithm to $w_1=(74, 125, 180)$ as before, and we again find 4 new
sets of invariants:
$$
( 174, 240, 246 ), ( 193, 85, 15 ), ( 268, 256, 143 ), ( 75, 263, 182 ). 
$$
In the end we expand the Igusa polynomials
\scriptsize
\begin{align*}
P_K &=
    X^{10} + 92X^9 + 72X^8 + 217X^7 + 98X^6 + 195X^5 + 233X^4 + 140X^3 +
        45X^2 + 123X + 171 \,,\\
Q_K &=
    X^{10} + 232X^9 + 195X^8 + 45X^7 + 7X^6 + 195X^5 + 173X^4 + 16X^3 +
        33X^2 + 247X + 237 \,,\\
R_K &=
    X^{10} + 240X^9 + 57X^8 + 213X^7 + 145X^6 + 130X^5 + 243X^4 + 249X^3 +
        181X^2 + 134X + 81 
\end{align*}
\normalsize
modulo $p=271$.
\end{example}

\normalsize
\begin{example} \label{sec7:ex2}
In the previous example, all the prime ideals of $K$ lying over $3$ gave rise to
an isogenous abelian surface. This phenomenon does not always occur. Indeed,
let $K$ be a primitive quartic CM-field and let $\gp_1,\ldots,\gp_n$ be the 
prime ideals of norm $3$. If we have a principally polarized
abelian surface $A/\Fp$ with endomorphism ring $\Order_K$, then the 
number of $(3,3)$-isogenous abelian surfaces with the same endomorphism ring
equals the cardinality of
$$
\{ m(\gp_1),\ldots,m(\gp_n)\}.
$$
There are examples where this set has {\it less\/} than $n$ elements. 

Take the cyclic field $K=\Q[X]/(X^4+219X^2+10512)$. The class group of $K$
is isomorphic to $\Z/2\Z \times \Z/2\Z$. The prime $3$ ramifies in $K$, and
we have $(3) = \gp_1^2\gp_2^2$. The primes $\gp_1,\gp_2$ in fact generate 
$\Cl(\Order_K)$. It is easy to see that for this field we have
$$
m(\gp_1) = m(\gp_2) \in \mathfrak{C}(K),
$$
so we only find {\it one\/} isogenous surface.
\end{example}

\begin{example} \label{sec7:bigCgroup}
Our algorithm is not restricted to cyclic CM-fields. In this example we 
let $K = \Q[X]/(X^4+22X^2+73)$ be a CM-field with Galois group~$D_4$. There
are $2$ equivalence classes of CM-types. We fix a CM-type $\Phi: K \rightarrow
\C^2$ and let $K_\Phi$ be the reflex field for $\Phi$. We have $K_\Phi =
\Q[X]/(X^4+11X^2+12)$, and $K$ and $K_\Phi$ have the same Galois closure~$L$.

As the real quadratic subfield $K^+ = \Q(\sqrt{3})$ has narrow class group
$\Z/2\Z$, the group $\mathfrak{C}(K)$ fits in an exact sequence
$$
1 \mapright{} \Z/2\Z \mapright{} \mathfrak{C}(K) \mapright{} \Z/4\Z
\mapright{} \Z/2\Z \mapright{} 1
$$
and a close inspection yields $\mathfrak{C}(K) \cong \Z/4\Z$. The prime $3$
factors as
$$
(3) = \gp_1 \gp_2 \gp_3^2
$$
in the reflex field, and we have $\Cl(\Order_{K_\Phi}) = \Z/4\Z = \langle
[\gp_1] \rangle$. The element $m(\gp_1) \in \mathfrak{C}(K)$ has order $4$,
and under the map $\mathfrak{C}(K) \mapright{f} \Cl(\Order_K) = \Z/4\Z$ the 
element $f(m(\gp_1))$ has order 2. We see that even though the ideal
$N_{L/K}(\gp_1\Order_L)$ has order 2 in the class group, the typenorm of $\gp_1$
has order~4. 

Of the 4 ideal classes of $K$, only 2 ideal classes are principally 
polarizable for $\Phi$. The other 2 ideal classes are principally polarizable
for `the other' CM-type. Furthermore, the two principally polarizable 
ideal classes each have {\it two\/} principal polarizations. 

The prime $p=1609$ splits completely in the Hilbert class field of $K_\Phi$. 
As in Example~\ref{sec7:ex1}, we do a random search to find that a surface 
$A/\Fp$ with Igusa invariants $w_0=(1563,789,704)\in\Fp^3$ has endomorphism ring
$\Order_K$. We apply Algorithm \ref{sec6:cm} to this point. As output, we get
$w_0$ again and two new points $w_1=(1396,1200,1520)$ 
and $w_2=(1350,1316,1483)$. The fact that we find $w_0$ again should come
as no surprise since $m(\gp_3) \in \mathfrak{C}(K)$ is the trivial element.
The points $w_1$ and $w_2$ correspond to $\gp_1$ and $\gp_2$. 

As expected we compute that the cycle
$$
w_0 = (1563,789,704) \mapright{\gp_1} (1396,1200,1520) \mapright{\gp_1}
(1276,1484,7) 
$$
$$
\mapright{\gp_1} (1350,1316,1483) \mapright{\gp_1} w_0
$$
has length~4. To find the full Igusa class polynomials modulo $p$, we do
a second random search. The remaining 4 points are 
$( 782, 1220, 257 )$, $( 1101, 490, 1321 )$, $( 577, 35, 471 )$, $( 1154, 723, 1456 )$.

\end{example}

\section{Obstruction to isogeny volcanos} \label{sec8}
For an ordinary elliptic curve $E/\Fp$ over a finite field, Kohel introduced~\cite{KohelThesis} an
algorithm to compute the endomorphism ring $\End(E)$, which has 
recently been improved in~\cite{BS}. One first computes
the endomorphism {\it algebra\/} $K$ by computing the trace of the Frobenius
morphism $\pi$ of $E$. If the index $[\Order_K:\Z[\pi]]$ is only divisible
by small primes $l$, then Kohel's algorithm uses the $l$-isogeny graph to determine the endomorphism ring.  The algorithm depends on the fact that the graph of $l$-isogenies looks like a `volcano', and one can quotient by subgroups of order $l$ to
move down the volcano until one hits the bottom. We refer to~\cite{FM, KohelThesis} for the details of this algorithm. This approach succeeds because of the following fact.

\begin{lemma} \label{sec8:lemma}
Let $E,E'/\Fp$ be two ordinary elliptic curves whose endomorphism rings are
isomorphic to the same order $\Order$ in an imaginary quadratic field~$K$. 
Suppose that $l\not =p$ is a prime such that the index $[\Order_K:\Order]$ is
divisible by~$l$. Then there are no isogenies of degree $l$ between $E$ 
and $E'$.
\end{lemma}
\Proof
This result is well known. Since the proof helps us understand what goes
wrong in dimension~2, we give the short proof. Suppose
that there does exist an isogeny $\varphi: E \rightarrow E'$ of degree~$l$.
By the Deuring lifting theorem~\cite[Thm.~13.14]{LaEF}, we can
lift $\varphi$ to an isogeny $\widetilde\varphi: \widetilde E \rightarrow
\widetilde E'$ defined over the ring class field for $\Order$. By CM-theory, we
can write $\widetilde E' = \C/I$ with $I$ an {\it invertible\/} $\Order$-ideal of norm $l$.
But since $l$ divides the index $[\Order_K : \Order]$, there are no invertible ideals
of norm~$l$.
\endproof\par
\ \par\noindent
Unlike the elliptic curve case, there are a greater number of possibilities for the endomorphism ring of an $(l,l)$-isogenous abelian surface $A/\Fp$.  
Necessarily, the order must contain $\ZZ[\pi,\overline{\pi}]$ where $\pi$ 
corresponds to the Frobenius endomorphism of~$A$.
Let $\varphi: A \rightarrow A'$ be an $(l,l)$-isogeny of principally polarized abelian surfaces where $\Order=\End(A)$ contains $\Order'=\End(A')$.  Since $\varphi$ splits multiplication by $l$, it follows that  $\ZZ+l\Order\subseteq\Order'\subseteq\Order$ and hence $\Order'$ has index dividing $l^3$ in $\Order$.  In addition, since the $\ZZ$-rank is greater than two, it is possible to have several non-isomorphic suborders of $\Order$ having the same index.

A natural question is whether the `volcano approach' for elliptic 
curves can be 
generalized to ordinary principally polarized abelian surfaces $A/\Fp$. The 
extension of Schoof's algorithm~\cite{Pila} enables us to compute 
the endomorphism algebra $K = \End(A)\otimes_\Z\Q$, and the problem is to compute the 
subring $\End(A) = \Order \subseteq \Order_K$. By working with explicit $l$-torsion
points for primes $l \mid  [\Order_K:\Z[\pi, \overline{\pi}]]$ one can determine
this subring, see~\cite{EisentragerLauter, FreemanLauter}. This 
approach requires working over large extension field of $\Fp$ and a
 natural question is whether we can generalize the volcano algorithm 
directly by using $(l,l)$-isogenies between abelian surfaces. However, the
analogous statement to Lemma~8.1 that \par
\medskip\noindent
\centerline{\it there are no $(l,l)$-isogenies between $A$ and $A'$}\par
\medskip\noindent
if $\End(A)$ and $\End(A')$ have isomorphic endomorphisms rings whose
conductor in $\Order_K$ divides $l$ does {\it not\/} hold in general. This is a theoretical obstruction to a straightforward generalization of the algorithm for elliptic curves.

We first give an example where the analogue of Lemma~\ref{sec8:lemma} for abelian surfaces fails.
\begin{example} \label{sec8:cycles}
Take the point $( 782, 1220, 257 )\in \F_{1609}^3$ which we found in Example 
\ref{sec7:bigCgroup}.  Below we depict the connected component of the
$(3,3)$-isogeny graph. The white dots represent surfaces with endomorphism
ring $\Order_K$, the black dots correspond to surfaces whose endomorphism ring is non-maximal.
The lattice of suborders of $\Order_K$ of $3$-power index that contain 
$\Z[\pi,\overline{\pi}]$ is completely described by the indices of the 
suborders in this case.
We have 
$\Z[\pi,\overline{\pi}] \subset \Order_{27}\subset \Order_9 \subset \Order_3 \subset \Order_K$, where the subscript 
denotes the index in~$\Order_K$.
\begin{center}
\includegraphics{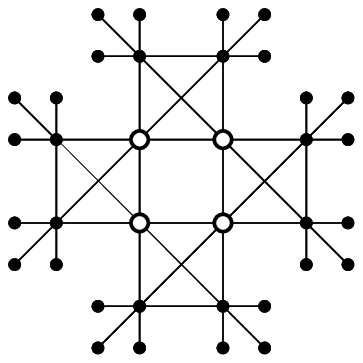}
\end{center}
The leaf nodes all have endomorphism ring $\Order_{27}$ and the remaining eight 
black vertices have endomorphism ring $\Order_3$. We observe that there are cycles in this graph other than at the `surface' of the volcano.\hfill$\diamond$
\end{example}
The reason that cycles can occur is the following.
Just like in
Lemma~\ref{sec8:lemma}, we can lift an isogeny $\varphi: A \rightarrow A'$ to
characteristic zero. By CM-theory, we can now write 
$\widetilde A = \C^2/\Phi(I)$ for
some invertible $\Order$-ideal~$I$. The isogenous surface $\widetilde A'$ equals
$\C^2 / \Phi(\ga^{-1}I)$ for an invertible $\Order$-ideal $\ga$ of norm~$l^2$. The
difference with the elliptic curve case is that there {\it do\/} exist 
invertible $\Order$-ideals of norm~$l^2$. Hence, the isogeny graph for 
abelian surfaces need not look like a `volcano'.

Another ingredient of the endomorphism ring algorithm for elliptic curves 
 can fail. In the
elliptic curve case, the following property of the $l$-isogeny graph is
essential.  Suppose that 
$E/\Fp$ has endomorphism ring $\Order$ and let $\varphi: E \rightarrow E'$ be
an isogeny from $E$ to an elliptic curve with endomorphism ring of index~$l$.
If $\varphi$ is defined over $\Fp$, then {\it all\/} $l+1$ isogenies of
degree $l$ are defined over $\Fp$.

The analogous statement for dimension 2 is that {\it all\/} $(l,l)$-isogenies
are defined over $\Fp$ as soon as there is one $(l,l)$-isogeny $\varphi: A
\rightarrow A'$ that is defined over $\Fp$. Here, $A'$ is an abelian 
surface with endomorphism ring of index dividing~$l^3$. This statement 
is {\it not\/} true, as the following example shows.

\begin{example} \label{sec8:nonregular}
Consider the cyclic quartic CM-field $K=\QQ[X]/(X^4+12X^2+18)$ 
which has class number $2$.  The Igusa class polynomials have degree $2$ and 
over $\F_{127}$ we find the corresponding moduli points $w_0=(118, 71, 63)$, 
$w_1=(98, 82, 56)$.
The isogeny graph is not regular:
\begin{center}
\includegraphics{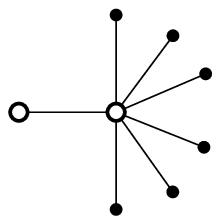}
\end{center}
The white dots represent the points having maximal endomorphism ring. There 
are $7$ points isogenous to $w_0$, which includes $w_1$. One cannot identify $w_1$ from the graph structure alone. This demonstrates that the isogeny graph is insufficient to determine the endomorphism rings; the polarized CM-lattices are also required.\hfill$\diamond$\par
\end{example}

The `reason' that the graph has this shape is the following. Let $\pi \in \Order_K$
correspond to the Frobenius morphism of a surface $A$ belonging to the vertex~$w_1$. If
$A'$ is $(l,l)$-isogenous to $A$, then $A'$ is defined over $\Fp$ if and only
if its endomorphism ring contains~$\pi$. Since there are several orders of
index dividing~$l^3$ in $\Order_K$, it can happen that $\pi$ is contained in some
of them, and not in others. In our example, the black points all have the
same endomorphism ring $\Order'$ with $\pi \in \Order'$. The 33 other isogenous surfaces
have an endomorphism ring that contains $\pi^3$, but not $\pi$. \par

\section*{Acknowledgements.} The authors thank Nils Bruin, Everett Howe, David Kohel, and John Voight for helpful dicussions and improvements to an early version the paper, the anonymous referee both for various detailed comments
and for motivating us to include a run time analysis, and Mike Rosen for 
proving~\cite{MR} a genus theory result that allowed us to prove Lemma~\ref{sec6:lemgenus}.


\begin{thebibliography}{99}

\bibitem{Bach}
E. Bach, \textit{Explicit bounds for primality testing and related problems}, Math. Comp., \textbf{55} (1990), 355–-380

\bibitem{BB}
W. L. Baily, A. Borel, \textit{Compactification of arithmetic quotients of bounded symmetric domains}, Ann. of Math., \textbf{84} (1966), 442--538.

\bibitem{BBEL}
J. Belding, R. Br\"oker, A. Enge, K. Lauter, \textit{Computing Hilbert class polynomials}, ANTS VIII, Springer LNCS \textbf{5011} (2008), 282--295.

\bibitem{BS}
G. Bisson, A. V. Sutherland, \textit{Computing the endomorphism ring of an ordinary elliptic curve over a finite field}. Preprint available at {\tt http://arxiv.org/abs/0902.4670}.

\bibitem{Magma}
W. Bosma, J. Cannon, C. Playoust, \textit{The Magma algebra system. I. The user language}, J. Symbolic Comput., {\bf 24} (1997), pp. 235--265.

\bibitem{BostMestre}
J.-B. Bost, J.-F. Mestre, \textit{Moyenne arithm\'etique-g\'eometrique et p\'eriodes de courbes de genre 1 et 2}, Gaz. Math. Soc. France, \textbf{38} (1988), 36--64.

\bibitem{BrokerLauterPrePub}
R. Br\"oker, K. Lauter, \textit{Modular polynomials for genus $2$}, LMS Jour. of Math. and Comp., \textbf{12} (2009), 326--339. 

\bibitem{CardonaQuer}
G. Cardona, J. Quer, \textit{Field of moduli and field of definition for curves of genus 2}, in ``Computational aspects of algebraic curves'', Lecture Notes Ser. Comput. 13, World Sci. Publ., Hackensack, NJ, 2005, 71-83.

\bibitem{3adic}
R. Carls,  D. Kohel, D. Lubicz, \textit{Higher-dimensional 3-adic CM construction},  J. Algebra, \textbf{319}  (2008), 971--1006.

\bibitem{CN}
C.-L. Chai, P. Norman, \textit{Bad reduction of the Siegel moduli scheme of genus two with $\Gamma_0(p)$-level structure}, Amer. J. Math., \textbf{122} (1990), 1003--1071.

\bibitem{DupontThesis}
R. Dupont, \textit{ Moyenne arithm\'etico-g\'eom\'etrique, suites de Borchardt et applications}, PhD thesis, \'Ecole polytechnique, 2006.


\bibitem{EisentragerLauter}
K. Eisentr\"ager, K. Lauter, \textit{A CRT algorithm for constructing genus 2 curves over finite fields}, Proceedings of Arithmetic, Geometry, and Coding Theory, (AGCT-10), 161--176.

\bibitem{FM}
M. Fouquet, F. Morain. \textit{Isogeny volcanoes and the SEA algorithm}, in Algorithmic Number Theory (2002), ANTS V, Springer LNCS {\bf 2369} (2002), pp. 276--291. 

\bibitem{FreemanLauter}
D. Freeman, K. Lauter, \textit{Computing endomorphism rings of Jacobians of genus 2 curves over finite fields}, in Algebraic Geometry and its Applications, World Scientific (2008), 29--66.

\bibitem{2adicFiveAuthor}
P. Gaudry, T. Houtmann, D. Kohel, C. Ritzenthaler, and A. Weng, \textit{The 2-adic CM method for genus 2 curves with application to cryptography}, Asiacrypt 2006 (Shanghai) Lect. Notes in Comp. Sci., 4284, 114--129, Springer-Verlag, 2006.

\bibitem{vdG82}
G. van der Geer, \textit{On the geometry of a {S}iegel modular threefold},
Math. Ann., \textbf{260} (1982), 317--350.

\bibitem{Gor97}
\hyphenation{Ma-nu-scrip-ta}
E. Z. Goren, \textit{On certain reduction problems concerning abelian surfaces},
Manuscripta Math., \textbf{94} (1997), no. 1, 33--43.

\bibitem{GorenLauter}
E. Z. Goren, K. Lauter, \textit{Genus 2 curves with complex multiplication}, preprint, available at {\tt http://arxiv.org/abs/1003.4759/}.

\bibitem{GruenewaldArticle}
D. Gruenewald, \textit{Computing Humbert surfaces and applications}, in  Arithmetic, geometry, cryptography and coding theory, Contemp. Math., \textbf{521}, Amer. Math. Soc., 2010, 59--69. 

\bibitem{Gruenewald}
D. Gruenewald, \textit{Explicit Algorithms for Humbert Surfaces}, PhD Thesis, University of Sydney, December, 2008.

\bibitem{IgusaAVM}
J. Igusa, \textit{Arithmetic variety of moduli for genus two},
Ann. of Math. (2), \textbf{72} (1960), 612--649.

\bibitem{IgusaMPI}
J. Igusa, \textit{Modular forms and projective invariants},
Amer. J. Math., \textbf{89} (1967), 817--855.

\bibitem{IgusaThetaConsts1}
J. Igusa, \textit{On the graded ring of theta-constants},
Amer. J. Math., \textbf{86} (1964), 219--246.

\bibitem{KM}
N. M. Katz, B. Mazur, \textit{Arithmetic moduli of elliptic curves}, Annals of Mathematics Studies 108, Princeton University Press (1985).

\bibitem{KohelThesis}
D. Kohel, \textit{Endomorphism rings of elliptic curves over finite fields}, PhD thesis, University of California, Berkeley, 1996.

\bibitem{Kohel}
D. Kohel, \textit{Complex multiplication and canonical lifts}, in Algebraic Geometry and its Applications, World Scientific (2008), 67--83.

\bibitem{LO}
J.~C. Lagarias and A.~M. Odlyzko, \textit{Effective versions of the {C}hebotarev
  density theorem}, Algebraic number fields: {$L$}-functions and {G}alois
  properties (Proc. Sympos., Univ. Durham, Duram, 1975), Academic Press, 1977,
  pp.~409--464.

\bibitem{LaEF}
S. Lang, \textit{Elliptic functions}, Springer GTM \textbf{112} (1987).

\bibitem{LaCM}
S. Lang, \textit{Complex multiplication}, Springer GTM \textbf{255} (1983).

\bibitem{Louboutin}
S. Louboutin, \textit{Explicit lower bounds for residues at $s=1$ of Dedekind
zeta functions and relative class numbers of CM-fields}, Trans. Amer. Math. Soc., \textbf{355} (2003), pp. 3079--3098.

\bibitem{Mestre}
J.-F. Mestre, \textit{Construction de courbes de genre 2 \`a partir de leurs modules}, in Effective methods in algebraic geometry, Birkh\"auser Progr. Math. \textbf{94} (1991), 313--334.

\bibitem{Mum}
D. Mumford, \textit{Abelian varieties}, Oxford University Press (1970).

\bibitem{JN}
J. Nakagawa, \textit{Orders of a quartic field}, Memoirs Amer. Math. Soc. \textbf{583} (1996).

\bibitem{MR}
M. Rosen, \textit{The $p$-rank of the class group in cyclic $p$-power extensions}, in preparation (2011).
 
\bibitem{Rung93}
B. Runge, \textit{On Siegel modular forms, part I}, J. Reine angew. Math., \textbf{436} (1993), pp. 57--85 

\bibitem{Pila}
J. Pila, \textit{Frobenius maps of abelian varieties and finding roots of unity in finite fields}, Math. Comp. {\bf 55} (1990), no. 192, pp. 745--763.

\bibitem{PoorYuen}
C. Poor and D. Yuen, \textit{Linear dependence among Siegel modular forms}, Math. Ann. 318 (2000), no. 2, 205--234

\bibitem{Shim}
G. Shimura, \textit{Abelian varieties with complex multiplication and
  modular functions}, Princeton Mathematical Series \textbf{46},
Princeton University Press New Jersey, 1998.


\bibitem{Siegel}
C.-L. Siegel, \textit{Symplectic Geometry}, Academic Press, New York, 1964.

\bibitem{Sil}
J. H. Silverman, \textit{Advanced topics in the arithmetic of elliptic curves}, Springer GTM \textbf{151} (1994).

\bibitem{Spallek}
A.-M. Spallek, \textit{Kurven vom {G}eschlecht 2 und ihre {A}nwendung in {P}ublic-{K}ey-{K}ryptosystemen},
PhD thesis, Universit\"at Gesamthochschule Essen, 1994.

\bibitem{psh}
P. Stevenhagen, \textit{The arithmetic of number rings}, Algorithmic Number Theory, Mathematical Sciences Research Institute Publications {\bf 44}, Cambridge University Press (2008), pp. 209--266.

\bibitem{Streng}
M. Streng, \textit{Complex multiplication of abelian surfaces}, PhD-thesis, Universiteit Leiden, 2010.

\bibitem{van Wamelen}
P. van Wamelen.
\textit{Examples of genus two {CM} curves defined over the rationals},
Math. Comp., 68 (1999), 307--320.

\bibitem{Weng}
A. Weng, \textit{Constructing hyperelliptic curves of genus 2 suitable for cryptography}, Math. Comp. 72 (2003), 435--458.


\end{thebibliography}
\end{document}